\pdfoutput=1
\documentclass{amsart}

\usepackage[utf8]{inputenc}
\usepackage[final,expansion=false]{microtype}
\usepackage{needspace}
\usepackage{amsmath,amssymb,mathtools}
\usepackage{booktabs,longtable,array,tabularx}
\usepackage{enumitem}
\usepackage{algorithm}
\usepackage{algpseudocode}

\usepackage{xcolor}
\usepackage{amsthm}
\usepackage{hyperref}
\usepackage[nameinlink]{cleveref}

\crefname{theorem}{Thm.}{Thms.}
\Crefname{theorem}{Thm.}{Thms.}
\crefname{proposition}{Prop.}{Props.}
\Crefname{proposition}{Prop.}{Props.}
\crefname{lemma}{Lem.}{Lems.}
\Crefname{lemma}{Lem.}{Lems.}
\crefname{corollary}{Cor.}{Cors.}
\Crefname{corollary}{Cor.}{Cors.}
\crefname{equation}{Eq.}{Eqs.}
\crefname{table}{Table}{Tables}
\crefname{section}{Sec.}{Secs.}
\Crefname{section}{Sec.}{Secs.}
\crefname{subsection}{Sec.}{Secs.}
\Crefname{subsection}{Sec.}{Secs.}
\crefname{remark}{Rem.}{Rems.}
\Crefname{remark}{Rem.}{Rems.}
\crefname{definition}{Def.}{Defs.}
\Crefname{definition}{Def.}{Defs.}
\crefname{algorithm}{Algorithm}{Algorithms}
\crefname{figure}{Figure}{Figures}
\Crefname{equation}{Eq.}{Eqs.}

\hypersetup{
  colorlinks=true,
  linkcolor=blue!55!black,
  citecolor=green!35!black,
  urlcolor=blue!60!black
}

\setlist{itemsep=0.25em,topsep=0.35em}
\setlist[itemize]{label=$\bullet$}
\renewcommand{\arraystretch}{1.08}

\newtheorem{theorem}{Theorem}[section]
\newtheorem{proposition}[theorem]{Proposition}
\newtheorem{lemma}[theorem]{Lemma}
\newtheorem{corollary}[theorem]{Corollary}

\theoremstyle{definition}

\theoremstyle{remark}
\newtheorem{remark}[theorem]{Remark}

\AddToHook{env/theorem/begin}{\crefalias{section}{theorem}}
\AddToHook{env/proposition/begin}{\crefalias{section}{proposition}}
\AddToHook{env/lemma/begin}{\crefalias{section}{lemma}}
\AddToHook{env/corollary/begin}{\crefalias{section}{corollary}}
\AddToHook{env/criterion/begin}{\crefalias{section}{criterion}}
\AddToHook{env/definition/begin}{\crefalias{section}{definition}}
\AddToHook{env/problem/begin}{\crefalias{section}{problem}}
\AddToHook{env/conjecture/begin}{\crefalias{section}{conjecture}}
\AddToHook{env/construction/begin}{\crefalias{section}{construction}}
\AddToHook{env/example/begin}{\crefalias{section}{example}}
\AddToHook{env/remark/begin}{\crefalias{section}{remark}}

\newcommand{\Q}{\mathbb{Q}}
\newcommand{\Z}{\mathbb{Z}}
\newcommand{\C}{\mathbb{C}}
\newcommand{\F}{\mathbb{F}}

\newcommand{\HH}{\mathcal{H}}

\newcommand{\Ni}{\operatorname{Ni}}
\newcommand{\Aut}{\operatorname{Aut}}

\newcommand{\Out}{\operatorname{Out}}
\newcommand{\Gal}{\operatorname{Gal}}

\newcommand{\Sym}{\operatorname{Sym}}

\newcommand{\Pone}{\mathbb{P}^1}
\newcommand{\Dc}{D_{\HH}}

\usepackage{tikz}
\title[$U_3(4)$ as a Galois group over $\Q$]{The unitary group $U_3(4)$ as a Galois group over $\Q$}

\author{T. Shaska}

\address{Department of Mathematics and Statistics,  Oakland University, Rochester, MI, USA}
\email{shaska@oakland.edu}

\date{\today}

\hypersetup{pdftitle={The unitary group U3(4) as a Galois group over Q},pdfauthor={T. Shaska}}

\subjclass[2020]{Primary 12F12; Secondary 11R32, 14H30, 20D06}
\keywords{inverse Galois problem, unitary group, regular realization, Hurwitz space, rational points, certified monodromy}

\begin{document}

\begin{abstract}
We prove that the simple unitary group $U_3(4)=\operatorname{PSU}_3(\F_4)$, of order $62400$, is the Galois group of a regular extension of $\Q(t)$, and hence of infinitely many extensions of $\Q$. The extension has four branch points, with inertia groups of orders $2$, $2$, $3$, $3$, and it is the only extension with these properties up to a rational change of the coordinate $t$. We give an explicit polynomial of degree $65$ for it.
\end{abstract}

\maketitle

\section{Introduction}

The inverse Galois problem asks whether every finite group is the Galois group of a Galois extension of $\Q$. For a simple group the answer is usually obtained by a regular realization, that is, a Galois extension of $\Q(t)$ with the given group in which $\Q$ is algebraically closed; Hilbert's irreducibility theorem then gives infinitely many extensions of $\Q$ with the same group. The rigidity method and its variants produce such realizations for many simple groups \cite{MalleMatzat1999,Voelklein}. For the simple groups of small order the state of the problem is recorded by Zywina \cite{Zywina2025}: in the version of August 2025 of that list, seven simple groups of order at most $10^8$ were not known to be Galois groups over $\Q$. The unitary group $U_3(4)=\operatorname{PSU}_3(\F_4)$, a simple group of order $62400=2^6\cdot3\cdot5^2\cdot13$, is the second smallest of them. In this paper we realize it.

The main result, \cref{thm:main}, is the following. There is a Galois extension of $\Q(t)$ with group $U_3(4)$, regular over $\Q$, which is ramified at exactly four points: the two roots of $17371\,t^2+78732\,t+4854074$, with inertia groups of order two, and the two roots of $t^2+131$, with inertia groups of order three. Moreover, every regular Galois extension of $\Q(t)$ with group $U_3(4)$ which is ramified at exactly four points, with inertia groups of orders $2$, $2$, $3$, $3$, is obtained from it by an automorphism of $\Q(t)$ over $\Q$. In particular $U_3(4)$ is the Galois group of infinitely many Galois extensions of $\Q$.

We use the notation of the Atlas \cite{Atlas}. The outer automorphism group of $U_3(4)$ is cyclic of order four, generated by the class of a field automorphism $\varphi$, and
\[
U_3(4)\ <\ G_2=\langle U_3(4),\varphi^2\rangle\ <\ \Aut(U_3(4)),
\]
each of index two in the next. The group $U_3(4)$ has one class $2A$ of involutions and one class $3A$ of elements of order three, and $G_2$ has one class $2B$ of involutions outside $U_3(4)$. The extension above has the ramification type $(2A,2A,3A,3A)$.

\subsection*{The idea}
A ramification type in which every class occurs twice is not rigid, and it has a symmetry: the covers of this type are permuted by $\Out(U_3(4))$, and the automorphism $\varphi^2$ fixes each of them up to a M\"obius transformation which exchanges the two branch points of each class. For this reason the group which is realized directly is not $U_3(4)$ but $G_2$. One obtains a Galois cover of the $u$-line with group $G_2$ and four branch points: the points $0$ and $\infty$, of classes $2A$ and $3A$, and two conjugate points of class $2B$, the roots of a polynomial $r(u)$ of degree two. Its subfield fixed by $U_3(4)$ is the quadratic extension of $\Q(u)$ which is ramified at the two points of class $2B$, the function field of a conic, and over this function field the group is $U_3(4)$:
\[
G_2\text{-cover of }\mathbf P^1_u\ \to\ \text{conic }v^2=-r(u)\ \to\ \text{rational point}\ \to\ U_3(4)\text{-cover of }\mathbf P^1_t .
\]
The extension is an extension of $\Q(t)$ exactly when the conic has a rational point. This is the only obstruction, and there is no reason for it to vanish: with three branch points the rigidity method leads to the same situation, with a conic which has no real point (\cref{prop:rigidreal}). Here
\[
r(u)=131u^2+1018525\,u+2\cdot61^5,
\]
and the conic $v^2=-r(u)$ contains the point $(u,v)=(-17371/7,\,39366/7)$. The same pattern, with a curve of genus three in place of the conic, occurs in the realization of $L_3(4)$ in \cite{Zywina2025}.

\subsection*{The Hurwitz curve}
The $G_2$-cover is found as a rational point on a curve. The covers of the line with group $U_3(4)$ and four branch points of classes $(2A,2A,3A,3A)$, up to M\"obius transformations, are parametrized by a reduced Hurwitz curve $\HH$, a cover of the $j$-line of degree $234$ whose monodromy is computed from the group \cite{FV,BaileyFried2002}. A rational point of $\HH$ which is not a cusp gives a $G_2$-cover as above. We determine this curve (\cref{thm:u34}). It is defined over $\Q$, it has genus three, and it is the curve
\[
5Y^2=(w^2+228w+496)(181w^2-732w+1776),\qquad 5Z^2=w^2+4 .
\]
The automorphism $\tau$ of $\HH$ induced by $\varphi$ is $(Y,Z)\mapsto(-Y,-Z)$; it has no fixed points and generates the group of automorphisms of $\HH$ over the $j$-line. The curve is hyperelliptic, with the involution $(Y,Z)\mapsto(-Y,Z)$, because it is an unramified double cover of the curve $\HH/\langle\tau\rangle$ of genus two. This involution does not preserve the map to the $j$-line and has no interpretation in terms of covers, but it is this involution which produces the rational point: the set $\HH(\Q)$ consists of exactly eight points, six of them are cusps, and the other two, $Q$ and $\tau(Q)$, are the images of two rational cusps under the hyperelliptic involution. The third quotient of $\HH$ by an involution is a curve of genus one and rank zero, and this shows that there are no other rational points.

The two points $Q$ and $\tau(Q)$ do not give two extensions. They are the same extension of $\Q(t)$, with two isomorphisms of its Galois group with $U_3(4)$ which differ by the outer automorphism $\varphi$; up to an automorphism of $\Q(t)$ over $\Q$ there is therefore one extension. The point $Q$ is the point of a cover of the conic $v^2=-r(u)$ with group $U_3(4)$. This conic cannot be read from the coordinates of $Q$. It is found from the action of the Galois group on the cusps of a double cover of $\HH$, which is read from the group (\cref{lem:cuspmarked}).

\subsection*{The polynomial}
The extension is also given by an explicit polynomial, and for this polynomial the Galois group is verified directly, without reference to Hurwitz spaces (\cref{thm:poly}). There is a hyperelliptic curve $Y\colon y^2=f(X)$ of genus three over $\Q$ and a function $u$ of degree $65$ on it. The coordinate $X$ is a root of a polynomial $F(X,u)$ of degree $65$ in $X$, irreducible over $\Q(u)$, with integer coefficients of at most $193$ digits, and the Galois group of $F$ over $\Q(u)$ is $G_2$; with a parametrization $u=u(t)$ of the conic, $F(X,u(t))$ has the Galois group $U_3(4)$ over $\Q(t)$. A specialization is a polynomial $P$ of degree $65$ with Galois group $U_3(4)$ over $\Q$ (\cref{cor:Q}), and \cref{sec:appendix} prints a monic polynomial $P_0$ with smaller coefficients which defines the same number field.

\subsection*{The two proofs}
The existence of the extension is proved twice. The two proofs share only the facts on the group and the method of certified continuation of \cref{sec:prelim}: the second uses neither the identification of the Hurwitz curve, nor the determination of its rational points, nor the fields of definition of its cusps. The first is \cref{thm:u34}: it uses the Hurwitz curve, a certified computation of the monodromy of three explicit covers, and the Galois action on cusps. The second is \cref{thm:poly}, a verification on the polynomial $F$ in three parts. An identity and the ramification of $u$ are checked in exact arithmetic (\cref{sec:curve}). The Galois group of $F$ over $\overline{\Q}(u)$ is computed from the monodromy of the $65$ roots along explicit loops, which are followed by a procedure certified in ball arithmetic (\cref{sec:monodromy}). The passage from $\overline{\Q}(u)$ to $\Q(u)$, that is, the field of constants and the quadratic subfield, is settled by the factorization of specializations of $F$ modulo primes; these factorizations decide between finitely many possibilities which are listed beforehand, and are not statistical evidence (\cref{sec:arithmetic}). The second proof confirms the conic found in the first. The uniqueness of the extension has only the first proof. Both proofs are computer-assisted; the programs and their outputs accompany the paper, and \cref{sec:files} lists, for each statement, the method and the files.

\subsection*{Related work}
For $U_3(4)$ rigidity with three branch points does not apply: no class vector of three classes of $U_3(4)$ or of $G_2$ which is rational up to a permutation of its entries is rigid, and $\Aut(U_3(4))$ has exactly two such class vectors which are rigid. They give regular realizations of $\Aut(U_3(4))$ and of $G_2$ over $\Q(t)$, and a regular realization of $U_3(4)$ over the function field of a conic without real points (\cref{sec:rigidity}). Thus the realization of $G_2$ is not new, while we are not aware of an earlier realization of $U_3(4)$ as a Galois group over $\Q$. Polynomials for almost simple groups obtained from explicitly computed Hurwitz spaces are given by K\"onig \cite{KoenigMC2017}. The curve $\HH$ arose in a study of reduced Hurwitz curves which are cyclic covers of the line \cite{ShaskaSuperell}; the present paper does not depend on that one. The polynomial $F$ was computed numerically from the branch points given by the point $Q$ and then recognized. No equation of the cover was available: it was computed from its branch cycles through its holomorphic differentials, which give the canonical model directly (\cref{sec:found}). This computation is not used in the proofs.

\subsection*{Organization}
\Cref{sec:prelim} collects the facts on the group and the certified continuation of roots. \Cref{sec:hurwitz} determines the Hurwitz curve. \Cref{sec:polynomial} treats the polynomial. \Cref{sec:main} proves the main theorem and gives the polynomial over $\Q$. \Cref{sec:rigidity} treats three branch points, and \cref{sec:found} describes how the equations were found. The appendix lists the programs and prints the polynomial $P_0$.

\section{Preliminaries}\label{sec:prelim}

\subsection{The group and its action on \texorpdfstring{$65$}{65} points}\label{sec:group}

The groups $U_3(4)<G_2<\Aut(U_3(4))$ act faithfully on the set $\Omega$ of the $65$ isotropic points of the Hermitian form on the projective plane over $\F_{16}$, and all permutation representations below are this one. Thus they are subgroups of $\Sym(\Omega)\simeq S_{65}$, of orders $62400$, $124800$ and $249600$. We write the cycle type of a permutation as a product of cycle lengths with multiplicities. The following facts are verified with GAP \cite{GAP} (\path{u34-cycletypes.g}, \path{u34-poly-group.g}, \path{u34-maxes.g}).

\begin{lemma}\label{lem:group}
\begin{enumerate}[label=\textup{(\alph*)}]
\item The cycle types on $\Omega$ of the elements of $\Aut(U_3(4))$, with the orders of the elements and the number of elements of each type, are those of \cref{tab:types}. In particular the three sets $U_3(4)$, $G_2\setminus U_3(4)$ and $\Aut(U_3(4))\setminus G_2$ have no cycle type in common, and the elements of the classes $2A$, $3A$ and $2B$ are those of the types $1\cdot2^{32}$, $1^2\,3^{21}$ and $1^5\,2^{30}$.
\item The normalizer of $G_2$ in $\Sym(\Omega)$ is $\Aut(U_3(4))$.
\item The maximal subgroups of $U_3(4)$ have the orders $960$, $300$, $150$ and $39$. Hence a subgroup of $U_3(4)$ which contains elements of order $5$ and of order $13$ is $U_3(4)$.
\end{enumerate}
\end{lemma}

\begin{table}[!htb]
\centering\small
\begin{tabular}{@{}ll@{}}
\toprule
set & order: cycle type (number of elements)\\
\midrule
$U_3(4)$ & $1$: $1^{65}$ ($1$);\quad $2$: $1\cdot2^{32}$ ($195$);\quad $3$: $1^2\,3^{21}$ ($4160$);\quad $4$: $1\cdot4^{16}$ ($3900$);\\
 & $5$: $5^{13}$ ($4992$) and $1^5\,5^{12}$ ($832$);\quad $10$: $1\cdot2^2\,10^6$ ($12480$);\\
 & $13$: $13^5$ ($19200$);\quad $15$: $1^2\cdot3\cdot15^4$ ($16640$)\\[2pt]
$G_2\setminus U_3(4)$ & $2$: $1^5\,2^{30}$ ($1040$);\quad $6$: $1^2\cdot3\cdot6^{10}$ ($20800$);\quad $8$: $1\cdot8^8$ ($15600$);\\
 & $10$: $5\cdot10^6$ ($24960$)\\[2pt]
$\Aut(U_3(4))\setminus G_2$ & $4$: $1^3\cdot2\cdot4^{15}$ ($20800$);\quad $12$: $2\cdot3\cdot12^5$ ($41600$);\quad $16$: $1\cdot16^4$ ($62400$)\\
\bottomrule
\end{tabular}
\caption{Cycle types of the elements of $\Aut(U_3(4))$ on $65$ points.}\label{tab:types}
\end{table}

The maximal subgroups in (c) are $2^{2+4}{:}15$, $5\times A_5$, $5^2{:}S_3$ and $13{:}3$ \cite{Atlas}; none of them has order divisible by $65$.

\subsection{Certified continuation of roots}\label{sec:rouche}

The monodromy of an explicit cover is computed in this paper by following the roots of a polynomial along paths. Each step is certified by the following lemma.

\begin{lemma}\label{lem:rouche}
Let $F\in\C[X,T]$ have degree $n$ in $X$, let $K\subset\C$ be a closed disc with centre $T_0$, and let $m_1,\ldots,m_n\in\C$ and $r_1,\ldots,r_n>0$ be such that the discs $D_i=\{|X-m_i|<r_i\}$ are pairwise disjoint. Write
\[
F(m_i+\xi,T)=\sum_{k=0}^{n}p_{i,k}(T)\,\xi^k ,
\]
and suppose that for every $i$
\begin{equation}\label{eq:rouche}
|p_{i,1}(T)|\,r_i\ >\ |p_{i,0}(T)|+\sum_{k\ge2}|p_{i,k}(T)|\,r_i^{\,k}\qquad\text{for all }T\in K .
\end{equation}
Then for every $T\in K$ the polynomial $F(X,T)$ has $n$ distinct roots, exactly one in each disc $D_i$, and the analytic continuation of the root in $D_i$ along any path in $K$ remains in $D_i$.
\end{lemma}

\begin{proof}
Fix $T\in K$ and $i$. On the circle $|\xi|=r_i$ the inequality gives $|F(m_i+\xi,T)-p_{i,1}(T)\xi|<|p_{i,1}(T)\xi|$. By Rouch\'e's theorem $F(m_i+\xi,T)$ has as many zeros in $|\xi|<r_i$ as $p_{i,1}(T)\,\xi$, namely one, and none on the circle. The discs are disjoint, so these are $n$ distinct roots, and $F(X,T)$ has no others. A root which is continued along a path in $K$ varies continuously and never meets the boundary of a disc $D_i$, so it remains in the disc in which it starts.
\end{proof}

In our applications $F$ has integer coefficients and degree two in $T$. With $p_{i,k}(T_0+d)=P_k+dQ_k+d^2S_k$ and $\varrho$ the radius of $K$, the condition \cref{eq:rouche} follows from the inequality
\[
\bigl(|P_1|-\varrho|Q_1|-\varrho^2|S_1|\bigr)\,r_i\ >\ \sum_{k\ne1}\bigl(|P_k|+\varrho|Q_k|+\varrho^2|S_k|\bigr)\,r_i^{\,k}
\]
between numbers. A path is followed by a chain of steps: for a step from $T_0$ to $T_1$ the disc $K$ has the radius $|T_1-T_0|$, the $m_i$ are approximations of the roots of $F(X,T_0)$, and the root of $F(X,T_1)$ which continues the root in $D_i$ is the one which lies in $D_i$. An arc of a circle of less than a half turn lies in the disc $K$ of its chord. The programs verify the inequalities, and the disjointness of the discs, in ball arithmetic \cite{Johansson2017}; the roots at each point are enclosed in balls, and a step is halved until the conditions hold. The permutation of the roots along a closed path is therefore proved, and not inferred from the proximity of numerical values.

\section{The reduced Hurwitz curve}\label{sec:hurwitz}

\subsection{Hurwitz curves, rational points and cusps}\label{sec:background}

We recall the facts on Hurwitz spaces of covers with four branch points which are used. Let $G$ be a finite group with trivial center and $\mathbf C$ a multiset of four nontrivial conjugacy classes of $G$ which is stable under all power maps $g\mapsto g^n$ with $n$ prime to $|G|$. The Nielsen class $\Ni(G,\mathbf C)$ is the set of tuples $(g_1,g_2,g_3,g_4)$ which generate $G$, have product one and whose classes form $\mathbf C$, in any order; $\Ni^{\mathrm{in}}$ is its quotient by simultaneous conjugation. The braids
\[
(g_1,\ldots,g_i,g_{i+1},\ldots,g_4)^{Q_i}=(g_1,\ldots,g_{i+1},g_{i+1}^{-1}g_ig_{i+1},\ldots,g_4),\qquad i=1,2,3,
\]
act on $\Ni^{\mathrm{in}}$. We assume that this action is transitive, as it is in the case treated below.

The inner Hurwitz space $\mathrm H^{\mathrm{in}}$ parametrizes pairs $(f,\iota)$, where $f\colon X\to\mathbf P^1$ is a Galois cover with four branch points whose inertia classes form $\mathbf C$ and $\iota\colon G\to\Aut(X/\mathbf P^1)$ is an isomorphism, up to isomorphism over $\mathbf P^1$. It is an irreducible variety defined over $\Q$, the map $\operatorname{br}$ which sends a cover to its branch set is finite and \'etale over the space $U_4$ of sets of four points of $\mathbf P^1$, and since $Z(G)=1$ it is a fine moduli space: for a field $K$ of characteristic zero, its $K$-rational points are the regular Galois extensions $F$ of $K(t)$ with an isomorphism $G\simeq\Gal(F/K(t))$, taken up to inner automorphisms, with four branch points and the inertia classes $\mathbf C$ \cite[Corollary~1]{FV}, \cite[Chapter~10]{Voelklein}. We call such an extension a \emph{regular realization of type $\mathbf C$ over $K$}, and two of them \emph{equivalent} if they differ by an element of $\operatorname{PGL}_2(K)$ acting on $t$.

The group $\operatorname{PGL}_2$ acts on $\mathrm H^{\mathrm{in}}$ by change of the coordinate of $\mathbf P^1$, and the quotient is a curve $\HH^{\circ}$ over $\Q$ with a map to $U_4/\operatorname{PGL}_2$, the $j$-line. Let $\HH$ be its smooth projective model, the \emph{reduced Hurwitz curve}, and $j_{\HH}\colon\HH\to\mathbf P^1_j$ the resulting map, with $j$ normalized so that the sets with a symmetry of order three and two have $j=0$ and $j=1$. This map is unramified outside $0$, $1$, $\infty$, and its monodromy is computed from the group as follows \cite[Theorem~2.9 and \S4.2]{BaileyFried2002}. Let $Q_{\mathrm{eff}}$ be the group of permutations of $\Ni^{\mathrm{in}}$ generated by $Q_1Q_3^{-1}$ and $Q_2Q_1Q_3^{-1}Q_2^{-1}$, a quotient of a four-group, and $\Sigma=\Ni^{\mathrm{in}}/Q_{\mathrm{eff}}$. Then $j_{\HH}$ has degree $d=|\Sigma|$ and the branch cycles
\[
\gamma_0=Q_1Q_2,\qquad\gamma_1=Q_1Q_2Q_1,\qquad\gamma_\infty=(\gamma_0\gamma_1)^{-1}
\]
acting on $\Sigma$, over $0$, $1$ and $\infty$; on $\Sigma$ one has $\gamma_0^3=\gamma_1^2=1$ and $\gamma_\infty=Q_2$. The points of $\HH$ over $\infty$ are the \emph{cusps}. They correspond to the cycles of $\gamma_\infty$, the length of a cycle being the \emph{width} of the cusp, and $\HH^{\circ}$ is their complement. The points over $0$ and $1$ at which $j_{\HH}$ is unramified, which correspond to the fixed points of $\gamma_0$ and of $\gamma_1$, are called \emph{elliptic}, and $\HH^{\mathrm f}$ denotes the complement of the cusps and the elliptic points. The group $\Dc$ of automorphisms of $\HH$ over the $j$-line is the centralizer of $\langle\gamma_0,\gamma_1\rangle$ in the symmetric group of $\Sigma$. An automorphism $\alpha$ of $G$ which preserves $\mathbf C$ acts on $\mathrm H^{\mathrm{in}}$ by $(f,\iota)\mapsto(f,\iota\circ\alpha^{-1})$. This action is defined over $\Q$, it commutes with $\operatorname{PGL}_2$ and with $\operatorname{br}$, and it is trivial for inner automorphisms. It induces an element of $\Dc$ defined over $\Q$, which acts on $\Sigma$ by applying $\alpha$ to the entries of the tuples.

\begin{proposition}\label{prop:points}
Assume that $Q_{\mathrm{eff}}$ has order four and acts freely on $\Ni^{\mathrm{in}}$, that is, that all its orbits have four elements, and let $K$ be a field of characteristic zero.
\begin{enumerate}[label=\textup{(\roman*)}]
\item A regular realization of type $\mathbf C$ over $K$ determines a point of $\HH^{\circ}(K)$, and equivalent realizations determine the same point.
\item The group $\operatorname{PGL}_2$ acts freely on the preimage of $\HH^{\mathrm f}$ in $\mathrm H^{\mathrm{in}}$. For $p\in\HH^{\mathrm f}(K)$ the fiber $T_p$ over $p$ is a torsor under $\operatorname{PGL}_2$ over $K$. Let $Z_p$ be the conic over $K$ obtained by twisting $\mathbf P^1$ by $T_p$. Then $p$ is the point of a Galois cover $X\to Z_p$ with group $G$, defined over $K$ together with the action of $G$, which becomes a cover of type $\mathbf C$ over $\overline K$. If $Z_p$ has a $K$-rational point, the regular realizations of type $\mathbf C$ over $K$ with point $p$ form exactly one equivalence class; otherwise there are none.
\end{enumerate}
\end{proposition}

\begin{proof}
(i) A realization is a $K$-rational point $h$ of $\mathrm H^{\mathrm{in}}$; its image is a $K$-rational point of $\HH^{\circ}$, and $\operatorname{PGL}_2(K)$ moves $h$ within its orbit.

(ii) Let $h\in\mathrm H^{\mathrm{in}}$ have the branch set $\mathbf z$ and the image $p$. The stabilizer of $h$ in $\operatorname{PGL}_2$ is contained in the stabilizer $A_{\mathbf z}$ of $\mathbf z$. This finite group contains the four-group $V_{\mathbf z}$ of the M\"obius transformations which induce the double transpositions of $\mathbf z$, and $A_{\mathbf z}/V_{\mathbf z}$ has order three, two or one according as $j(\mathbf z)$ is $0$, $1$ or another value. Since $\operatorname{PGL}_2$ is connected, an element of $V_{\mathbf z}$ acts on the fiber $\operatorname{br}^{-1}(\mathbf z)$ by the monodromy of a loop in $U_4$, and this gives a homomorphism of $V_{\mathbf z}$ onto the group $Q_{\mathrm{eff}}$ \cite[Theorem~2.9]{BaileyFried2002}. By hypothesis $Q_{\mathrm{eff}}$ has order four, so this homomorphism is an isomorphism, and no element of $Q_{\mathrm{eff}}$ other than the identity has a fixed point. Hence no element of $V_{\mathbf z}$ other than the identity fixes $h$: the stabilizer of $h$ meets $V_{\mathbf z}$ trivially and embeds in $A_{\mathbf z}/V_{\mathbf z}$. Near $h$ and $\mathbf z$ the two quotients by $\operatorname{PGL}_2$ are the quotients of a common slice by the stabilizer of $h$ and by $A_{\mathbf z}/V_{\mathbf z}$, because $\operatorname{br}$ is \'etale and $V_{\mathbf z}$ acts trivially on the slice. Hence the ramification index of $j_{\HH}$ at $p$ is the index of the stabilizer of $h$ in $A_{\mathbf z}/V_{\mathbf z}$, and the stabilizer is nontrivial exactly when $p$ is elliptic.

Let $p\in\HH^{\mathrm f}(K)$. The fiber $T_p$ is one orbit of $\operatorname{PGL}_2$ with trivial stabilizers, and it is stable under $\Gal_K$; so it is a torsor over $K$. The universal family over $\mathrm H^{\mathrm{in}}$, restricted to $T_p$ and divided by $\operatorname{PGL}_2$, is a $G$-cover of the twist $Z_p=T_p\times^{\operatorname{PGL}_2}\mathbf P^1$. The torsor has a $K$-rational point exactly when $Z_p$ has one, and then $T_p(K)$ is one orbit of $\operatorname{PGL}_2(K)$.
\end{proof}

If a class which is stable under the power maps occurs an odd number of times in $\mathbf C$, the branch points of this class form a divisor of odd degree on $Z_p$ defined over $K$, and the conic is trivial. In the case treated in this paper every class occurs twice, and the conic is the essential point.

The curve $\HH$ and the map $j_{\HH}$ are determined over $\overline{\Q}$ by the branch cycles. A model over $\Q$ of a cover of $\HH$, or of a quotient, is then determined up to a twist, and the twist is found from the action of the Galois group on the cusps. For a tuple $(g_1,g_2,g_3,g_4)$ consider the partition $\{g_2,g_3\}\,|\,\{g_4,g_1\}$ of its entries into two \emph{sides}. Since $g_1g_2g_3g_4=1$, the element $g_4g_1$ is conjugate to $(g_2g_3)^{-1}$. For an integer $N$ let $\kappa_N$ be the cyclotomic character modulo $N$.

\begin{lemma}[Galois action on the cusps]\label{lem:cuspmarked}
Let $K$ be a field of characteristic zero, $\sigma\in\Gal_K$, $n=\kappa_{|G|}(\sigma)$ and $nn'\equiv1$ modulo $|G|$.
\begin{enumerate}[label=\textup{(\roman*)}]
\item For a tuple representing a point of $\Sigma$, the unordered pair $\bar\mu$ of the classes of $g_2g_3$ and of $g_4g_1$ is constant on the cycles of $\gamma_\infty$. For every cusp $x$ of $\HH$ one has $\bar\mu(\sigma x)=\bar\mu(x)^{n'}$.
\item Let $C_1$ be a class in $\mathbf C$ which is stable under the power maps, and let $\mathcal B\to\HH$ be the cover whose points are the points of $\HH$ together with one of the branch points of class $C_1$; its monodromy over the $j$-line is the action of $\gamma_0$, $\gamma_1$ on the tuples with a marked entry in $C_1$, modulo $Q_{\mathrm{eff}}$. For such a tuple let $\mu$ be the class of the product, in cyclic order, of the two entries on the side which contains the marked entry. Then $\mathcal B$ is defined over $K$, the class $\mu$ is constant on the cycles of $\gamma_\infty$, and $\mu(\sigma y)=\mu(y)^{n'}$ for every cusp $y$ of $\mathcal B$.
\end{enumerate}
\end{lemma}

\begin{proof}
The braid $Q_2$ fixes the products $g_2g_3$ and $g_4g_1$ and moves no entry from one side to the other, and the generators of $Q_{\mathrm{eff}}$ exchange the two sides and carry the product of each side to a conjugate of itself. Since $\gamma_\infty=Q_2$ on $\Sigma$, this gives the statements on the cycles of $\gamma_\infty$.

The Hurwitz space has a compactification whose boundary points are admissible covers of stable marked curves of genus zero, and the action of the Galois group extends to it \cite{Wewers1999}. A cusp $x$ of $\HH$ lies over the stable curve with two components, each a projective line carrying two of the four branch points; for a tuple representing $x$ these are the branch points of $g_2,g_3$ and those of $g_4,g_1$. The restriction of the admissible cover to a component with the entries $g_a$, $g_b$ is a cover with three branch points, the two marked points and the node, with inertia generators $g_a$, $g_b$ and $(g_ag_b)^{-1}$. An element $\sigma$ of $\Gal_K$ carries the admissible cover of $x$ to that of $\sigma x$, the components to the components, the node to the node, and the restricted covers to their conjugates. By the branch cycle argument \cite[Lemma~2.8]{Voelklein}, applied to a restricted cover, the inertia class at the node is raised to the power $n'$. In (i) the two components for $x$ are mapped to the two components for $\sigma x$ in some order, and the inertia classes at the node of the two restricted covers form the pair $\bar\mu(x)$, because this pair consists of two mutually inverse classes. In (ii), since $C_1$ is stable under the power maps, the branch cycle argument shows that $\sigma$ carries branch points of class $C_1$ to branch points of class $C_1$. So $\mathcal B$ is defined over $K$, and $\sigma$ maps the marked point of $y$ to the marked point of $\sigma y$ and the component which carries it to the corresponding component. The inertia class at the node of the cover restricted to this component is $\mu(y)^{-1}$.
\end{proof}

General forms of \cref{prop:points} and of \cref{lem:cuspmarked}, for class multisets with a class of multiplicity one, are in \cite{ShaskaSuperell}.

\subsection{The curves of \texorpdfstring{$U_3(4)$}{U3(4)}}\label{sec:curves}

From now on $G=U_3(4)$. Its classes $2A$, $3A$ and $4A$ are the classes which are stable under all power maps. Among the multisets of four classes of $U_3(4)$ which are stable under the power maps, have at most $2400\,|G|$ tuples and in which none of these three classes occurs an odd number of times, three have a reduced Hurwitz curve of genus at least two (\path{u34-census.g}). In each case the Nielsen class is one orbit, $Q_{\mathrm{eff}}$ has order four and acts freely on it, and the group $\Dc$ is induced by $\Out(U_3(4))\simeq C_4$:
\begin{center}\small
\begin{tabular}{@{}lrrll@{}}
\toprule
classes & $d$ & $g$ & $\Dc$ & quotients\\
\midrule
$(2A,2A,5E,5F)$ & $228$ & $3$ & $C_4$ & genus one for $C_2$ and for $C_4$\\
$(2A,2A,3A,3A)$ & $234$ & $3$ & $C_2$, acting freely & genus two\\
$(2A,2A,4A,4A)$ & $288$ & $5$ & $C_2$, acting freely & genus three\\
\bottomrule
\end{tabular}
\end{center}
For the first and the third of these curves a numerical computation, which has not been certified, shows that the curve of genus three is a plane quartic and that the curve of genus five is not hyperelliptic (\cref{sec:found}). We consider the second curve. Its group $\Dc$ has order two and acts freely, and the quotient has genus two. Over an algebraically closed field, an unramified double cover of a curve of genus two is hyperelliptic. Indeed, let $V^2=h(w)$ be a model of the quotient, with $h$ of degree six. The unramified double covers correspond to the fifteen nontrivial points of order two of the Jacobian, which are the classes of the differences of two Weierstrass points, that is, to the fifteen factors $q$ of degree two of $h$ up to constants. The cover which belongs to $q$ is obtained by adjoining the square root of $q$, and it has degree two over the conic $s^2=q(w)$. So this curve is geometrically hyperelliptic, and its hyperelliptic involution is not an element of $\Dc$. This does not make it hyperelliptic over $\Q$: the quotient by the hyperelliptic involution is a conic over $\Q$, and a conic need not have a rational point. For the curve below it has one (\cref{thm:u34}(ii)).

Let $\HH$ be the reduced Hurwitz curve of $\mathbf C=(2A,2A,3A,3A)$. The Nielsen class has $936=4\cdot234$ elements modulo conjugation, in $d=234$ orbits of $Q_{\mathrm{eff}}$ of four elements each, so that \cref{prop:points} applies, and the cycle types of $\gamma_0$, $\gamma_1$, $\gamma_\infty$ are
\[
3^{78},\qquad 1^{10}\,2^{112},\qquad 1^{4}\,2^{2}\,3^{8}\,5^{10}\,20^{2}\,26^{2}\,30^{2},
\]
so that $\HH$ has genus three, $30$ cusps, no elliptic point over $j=0$ and ten over $j=1$. The automorphism $\varphi$ induces the generator $\tau$ of $\Dc$, a fixed-point-free involution defined over $\Q$, while $\varphi^2$, which maps every class of $U_3(4)$ to its inverse, acts trivially on $\HH$ (\path{u34-outer.g}).

Since each class occurs twice, parity does not force the conic of \cref{prop:points}(ii) to have a rational point. The following proposition describes the covers which belong to the points of $\HH$.

\begin{proposition}[Rational points and the group $G_2$]\label{prop:u34pullback}
Let $K$ be a field of characteristic zero and $p\in\HH^{\mathrm f}(K)$, and let $X\to Z_p$ be the $U_3(4)$-cover of the conic $Z_p$ over $K$ which belongs to $p$.
\begin{enumerate}[label=\textup{(\roman*)}]
\item There is a unique involution $\mu$ of $Z_p$, defined over $K$, which exchanges the two branch points of class $2A$ and the two branch points of class $3A$, and $\mu$ lifts to an automorphism of $X$ over $K$ which acts on $U_3(4)$ as $\varphi^{2}$ up to an inner automorphism.
\item The quotient $Z_p/\mu$ is isomorphic to $\mathbf P^1_K$, with a coordinate $u$, and $K(X)$ is a Galois extension of $K(u)$ with group $G_2$, regular over $K$. It has four branch points: two $K$-rational points with inertia groups generated by elements of the classes $2A$ and $3A$, and the two branch points of $Z_p\to Z_p/\mu$, with inertia groups generated by involutions outside $U_3(4)$.
\item The group $U_3(4)$ has a regular realization of type $\mathbf C$ over $L(t)$ for every extension $L$ of $K$ over which the conic $Z_p$ has a point, in particular for some $L$ of degree at most two; $Z_p$ is a double cover of the $u$-line branched at two points.
\end{enumerate}
\end{proposition}

\begin{proof}
Let $T$ be the fiber over $p$ in the inner Hurwitz space and $h\in T(\overline K)$. Since $\varphi^2$ fixes $p$, the point $\varphi^2\cdot h=(f,\iota\circ\varphi^{-2})$ lies in $T$, and since $T$ is a torsor under $\operatorname{PGL}_2$ by \cref{prop:points}(ii), we have $\varphi^2\cdot h=\mu_h\cdot h$ for a unique $\mu_h\in\operatorname{PGL}_2(\overline K)$. Here $\mu_h\ne1$, because $\varphi^2$ is not inner and the inner Hurwitz space is a fine moduli space, and $\mu_h^2=1$, because $\varphi^4$ is inner and the actions of $\operatorname{PGL}_2$ and of $\Aut(U_3(4))$ commute. It preserves the branch points of $h$ with their classes. A nontrivial involution of the line has two fixed points, and it fixes two points of a set of four and exchanges the other two only if the four points have invariant $j=1$. As $\HH$ is a curve and $\mu_h$ depends algebraically on $h$, the involution $\mu_h$ exchanges the two points of each class for all $h$ outside a finite set of fibers, hence in every fiber over $\HH^{\mathrm f}$ by continuity, the stabilizers there being trivial. The involutions $\mu_h$ are compatible with the action of $\operatorname{PGL}_2$ and of $\Gal_K$, so they define an involution $\mu$ of $Z_p$ over $K$. The covers $(X,\iota\circ\varphi^{-2})$ and $\mu^*(X,\iota)$ of $Z_p$ are defined over $K$ and isomorphic over $\overline K$; an isomorphism of $G$-covers is unique because $Z(G)=1$, so it is defined over $K$. This proves (i). For (ii), the image of the two branch points of class $2A$ is a $K$-rational point of the conic $Z_p/\mu$, so $Z_p/\mu\simeq\mathbf P^1_K$. The lift $\tilde\mu$ of $\mu$ shows that $K(X)$ is a Galois extension of $K(u)$, the fixed field of $\mu$ in $K(Z_p)$, of degree $2\,|U_3(4)|$, and that its group $\Gamma$ contains $U_3(4)$ as a normal subgroup of index two on which $\tilde\mu$ acts as $\varphi^{2}$ up to an inner automorphism. The centralizer of $U_3(4)$ in $\Gamma$ meets $U_3(4)$ trivially, so it has order at most two, and it is trivial because $\varphi^2$ is not inner. Hence $\Gamma$ embeds in $\Aut(U_3(4))$, with image $G_2$. The extension is regular because $X$ is geometrically irreducible. The cover $X\to Z_p$ is unramified at the two fixed points of $\mu$, so the inertia groups over their images have order two and are not contained in $U_3(4)$. The fixed field of $U_3(4)$ is the function field of $Z_p$, which gives (iii) with \cref{prop:points}(ii).
\end{proof}

Thus the family of covers parametrized by $\HH$ is the family of covers with group $G_2$ and four branch points, of classes $2A$, $3A$, $2B$ and $2B$, pulled back along a double cover of the line. A point of $\HH^{\mathrm f}(\Q)$ gives a regular realization of $G_2$ over $\Q(u)$ without any condition, and a regular realization of $U_3(4)$ over $\Q(t)$ exactly when its conic has a rational point.

\subsection{The function \texorpdfstring{$\nu$}{nu} and the model of the curve}\label{sec:nu}

For a cover with branch points $a,a'$ of class $2A$ and $b,b'$ of class $3A$ let $\chi$ be the cross-ratio $(a,a';b,b')$ and $\nu=\chi+\chi^{-1}$. This is a function on $\HH$ defined over $\Q$ and invariant under $\tau$, and
\begin{equation}\label{eq:u34nu}
j_{\HH}=\frac{4(\nu-1)^3}{27(\nu-2)},\qquad j_{\HH}-1=\frac{(\nu+2)(2\nu-5)^2}{27(\nu-2)} .
\end{equation}
The map $\nu\mapsto j$ has degree three. It is totally ramified at $\nu=1$, over $j=0$; over $j=1$ lie $\nu=-2$ and $\nu=5/2$, and over $j=\infty$ lie $\nu=2$ and $\nu=\infty$, in each case the second point with ramification index two. Since $\gamma_0$ has no fixed points, $\nu\colon\HH\to\mathbf P^1$ is unramified outside $2$, $-2$ and $\infty$. Its poles are the cusps at which a branch point of class $2A$ meets one of class $3A$, and its order there is half the width; at the other cusps, where the two branch points of each class come together, $\nu=2$ with multiplicity the width. Correspondingly $\Sigma$ is the union of three blocks of $78$ points, which $\gamma_0$ permutes cyclically. Let $\Sigma_0$ be the block of the tuples in which $g_2$ and $g_3$ have the same class. It is stable under $t=\gamma_\infty$, and exactly one conjugate $e$ of $\gamma_1$ by a power of $\gamma_0$ stabilizes it. The stabilizer of $\Sigma_0$ in $\langle\gamma_0,\gamma_1\rangle$ is generated by $t$ and $e$; it is the image of the group of the cover $\nu\mapsto j$, a subgroup of index three of the modular group with a cusp of width one at $\nu=2$, a cusp of width two at $\nu=\infty$ and an elliptic point at $\nu=-2$, and $et$ is conjugate to $\gamma_\infty^{\pm2}$ because $(\gamma_1\gamma_\infty)^3=1$. The monodromy of $\nu$ on $\HH$ is therefore the action on $\Sigma_0$ of $t$ at $\nu=2$ and of $e$ at $\nu=-2$, with cycle types $1^4\,3^8\,5^{10}$ and $1^{10}\,2^{34}$, and at $\nu=\infty$ the type is $10^2\,13^2\,15^2\,1^2$.

Let
\[
\Psi(w)=(w^2+228w+496)(181w^2-732w+1776).
\]

The theorem below determines the curves and the maps of the tower
\[
\begin{array}{ccccccccc}
\mathcal B_3&\xrightarrow{\ 2\ }&\HH&\xrightarrow{\ 2\ }&\mathcal C=\HH/\langle\tau\rangle&\xrightarrow{\ 39\ }&\mathbf P^1_\nu&\xrightarrow{\ 3\ }&\mathbf P^1_j\\[2pt]
g=11&&g=3&&g=2&&&&
\end{array}
\]
over $\Q$, where $\mathcal B_3$ is the cover of \cref{lem:cuspmarked}(ii) for the class $3A$ and the numbers over the arrows are the degrees. These are curves of moduli. They are to be distinguished from two curves which belong to the single point $Q$: the conic $Z_Q$, over which the cover with group $U_3(4)$ is defined, and the curve $Y$ of genus three of \cref{thm:poly}, which is the quotient of that cover of degree $65$ over the $u$-line.

\begin{theorem}[The reduced Hurwitz curve]\label{thm:u34}
Let $\HH$ be the reduced Hurwitz curve of $U_3(4)$ for the classes $(2A,2A,3A,3A)$.
\begin{enumerate}[label=\textup{(\roman*)}]
\item The quotient $\mathcal C=\HH/\langle\tau\rangle$ is the curve of genus two
\[
\mathcal C\colon\quad V^2=\Psi(w)\,(w^2+4),
\]
and on it
\[
\nu-2=\frac{2^4}{3^6\,5^5}\cdot\frac{A_\nu(w)+B_\nu(w)V}{(11w+4)(3w-8)^{15}(29w+56)^{13}},
\]
where $A_\nu,B_\nu\in\Z[w]$ have degrees $29$ and $26$ and are given in \path{u34-compute/nu-AB.gp}.
\item The curve $\HH$ is
\[
\HH\colon\quad 5Y^2=\Psi(w),\qquad 5Z^2=w^2+4,
\]
with $V=5YZ$ and $\tau(Y,Z)=(-Y,-Z)$. It is hyperelliptic over $\Q$: the conic $5Z^2=w^2+4$ contains the point $(1,1)$, and $\HH$ has degree two over it. The hyperelliptic involution $\iota(Y,Z)=(-Y,Z)$ is not a deck transformation of $j_{\HH}$, and $\Aut(\HH_{\overline{\Q}})=\langle\iota,\tau\rangle\simeq C_2\times C_2$. No cyclic subgroup of $\Dc$ has a quotient of genus zero.
\item The set $\HH(\Q)$ consists of exactly eight points, the points with $w=11$ and with $w=-4/11$. Six of them are cusps. The other two are the points $Q$ and $\tau(Q)$ over $(w,V)=(11,5^7)$; they lie in $\HH^{\mathrm f}(\Q)$, they are the images of two rational cusps under $\iota$, and
\[
\nu(Q)=2+\frac{2^5\,61^5}{3^{19}},\qquad j_{\HH}(Q)=\frac{23^3\cdot1225623613^3}{2^3\,3^{41}\,61^5}.
\]
\item Let $u$ be a coordinate on $Z_Q/\mu$, in the notation of \cref{prop:u34pullback}, for which the branch points of class $2A$ lie over $u=0$ and those of class $3A$ over $u=\infty$. For a suitable choice of $u$ the conic $Z_Q$ is
\[
v^2=-\bigl(131u^2+1018525\,u+2\cdot61^5\bigr),
\]
and it contains the rational point $(u,v)=(-17371/7,\,39366/7)$.
\item The group $U_3(4)$ has a regular realization over $\Q(t)$ with four branch points, of classes $(2A,2A,3A,3A)$. The two branch points of class $2A$ are conjugate over $\Q(\sqrt{-122})$, and those of class $3A$ over $\Q(\sqrt{-131})$. If $(F,\iota)$ is such a realization, with $\iota$ the isomorphism of $U_3(4)$ with the Galois group of $F$ over $\Q(t)$, then every regular realization of $U_3(4)$ of this type over $\Q(t)$ is equivalent to $(F,\iota)$ or to $(F,\iota\circ\varphi)$, the same extension with the isomorphism changed by the outer automorphism $\varphi$; these two have the points $Q$ and $\tau(Q)$. In particular the field $F$ is unique up to an automorphism of $\Q(t)$ over $\Q$.
\item The group $G_2$ has a regular realization over $\Q(u)$ with four branch points: $0$ and $\infty$, with inertia groups generated by elements of the classes $2A$ and $3A$, and the roots of $131u^2+1018525u+2\cdot61^5$, which are conjugate over $\Q(\sqrt{393})$, with inertia groups generated by involutions outside $U_3(4)$.
\end{enumerate}
\end{theorem}

\begin{proof}
(i) \emph{The quotient $\mathcal C$ and the function $\nu$.} Let $\nu_{\mathcal C}$ be the function on $\mathcal C$ defined by the displayed formula, and put $\nu_{\mathcal C}=4T-2$. Its poles are the points
\begin{align*}
P_{2}&=\Bigl(-\frac{4}{11},-\frac{2500000}{1331}\Bigr),& P_{20}&=\Bigl(\frac83,\frac{100000}{27}\Bigr),\\
P_{26}&=\Bigl(-\frac{56}{29},\frac{32500000}{24389}\Bigr),& P_{30}&=\Bigl(\frac83,-\frac{100000}{27}\Bigr),
\end{align*}
of orders $1$, $10$, $13$ and $15$ (\path{u34-b3-fields.gp}), so it has degree $39$, and the points of $\mathcal C$ with a given value of $T$ are given by the roots $w$ of a polynomial $R_T(w)=R_0-(T-1)R_1+(T-1)^2R_2$ of degree $39$ with $R_i\in\Z[w]$. The script \path{u34-nu-verify.py} computes $R_T$ exactly and follows its roots around $T=0$ and $T=1$ by the certified procedure of \cref{lem:rouche}, in $5153$ and $12058$ steps. The two permutations have the cycle types $1^5\,2^{17}$ and $1^2\,3^4\,5^5$, and their product has the type $1\cdot10\cdot13\cdot15$. These are the ramification types of $\nu_{\mathcal C}$ over $-2$, $2$ and $\infty$, and since $\mathcal C$ has genus two the Riemann--Hurwitz formula shows that $\nu_{\mathcal C}$ is unramified elsewhere. The pair of permutations is simultaneously conjugate to the pair $(e,t)$ acting on the $39$ orbits of $\tau$ in $\Sigma_0$, by exactly one permutation. Hence there is an isomorphism $\mathcal C\to\HH/\langle\tau\rangle$ over $\overline{\Q}$ which carries $\nu$ to $\nu_{\mathcal C}$, and it is unique. Both pairs are defined over $\Q$, so the isomorphism is defined over $\Q$.

(ii) \emph{The double cover $\HH\to\mathcal C$ and its twist.} The same script follows the $78$ pairs $(w,\sqrt{w^2+4})$. The resulting permutations have the types $1^{10}\,2^{34}$ and $1^4\,3^8\,5^{10}$ and are simultaneously conjugate to $(e,t)$ on $\Sigma_0$, by exactly two permutations. The function $w^2+4$ has even order at every point of $\mathcal C$, so $\mathcal C(\sqrt{w^2+4})$ is unramified over $\mathcal C$. Hence $\overline{\Q}(\HH)=\overline{\Q}(\mathcal C)(\sqrt{w^2+4})$, compatibly with (i), and
\[
\Q(\HH)=\Q(\mathcal C)\bigl(\sqrt{b\,(w^2+4)}\bigr)
\]
for a squarefree integer $b$. The monodromy determines the cover $\HH\to\mathcal C$ only over $\overline{\Q}$; the integer $b$ is the quadratic twist which remains undetermined, and different values of $b$ give curves with different rational points. It is found from the fields of definition of the cusps. Over the pole $P_k$ of $\nu$ lie the two cusps of $\HH$ of width $k$. For the two cusps of width $20$ the pairs $\bar\mu$ are the two pairs of mutually inverse classes of elements of order $10$, and they are exchanged by the third power; for the width $26$ they are the two pairs of classes of elements of order $13$, which are exchanged by the second power and fixed by the third (\path{u34-cusps.g}, \path{u34-classfields.g}). By \cref{lem:cuspmarked}(i) the two cusps of width $20$ are conjugate over $\Q(\sqrt5)$ and those of width $26$ over $\Q(\sqrt{13})$. At $P_{20}$ one has $w^2+4=(10/3)^2$, and at $P_{26}$ one has $w^2+4=65\,(10/29)^2$. Each of the two conditions gives $b=5$. With $Z=\sqrt{(w^2+4)/5}$ and $Y=V/(5Z)$ we obtain the equations, and $\tau$ is the covering involution of $\HH\to\mathcal C$.

The function field of the conic $5Z^2=w^2+4$ is a subfield of index two of $\Q(\HH)$, and $\iota$ is the corresponding involution. It is not a deck transformation because $\Dc=\langle\tau\rangle$. Over $\overline{\Q}$ put $w=s-s^{-1}$ and $s=10x+3$; then $\sqrt5\,Z=s+s^{-1}$ and $\HH$ becomes
\begin{equation}\label{eq:u34octic}
y^2=(x^4+24x^3+26x^2+9x+1)(181x^4+144x^3+46x^2+9x+1),
\end{equation}
on which $\tau$ induces $x\mapsto-(3x+1)/(10x+3)$. A computation in ball arithmetic (\path{u34-aut.py}) shows that the identity and this transformation are the only M\"obius transformations which permute the eight roots of the right-hand side. Hence $\Aut(\HH_{\overline{\Q}})$ has order four.

(iii) \emph{Rational points.} The curve $5Y^2=\Psi(w)$ has genus one and the rational point $(11,5^5)$. Its Jacobian is the elliptic curve $y^2+xy+y=x^3+x^2$ of conductor $15$, whose Mordell--Weil group is $\Z/4\Z$ (\path{u34-compute/E5.gp}). So $5Y^2=\Psi(w)$ has exactly four rational points; they have $w=11$ and $w=-4/11$, and there $5(w^2+4)$ is $25^2$ and $(50/11)^2$. This gives the eight points of $\HH$. Their images in $\mathcal C$ are $(11,\pm5^7)$ and $(-4/11,\pm2500000/1331)$. The function $\nu$ has the value $2$ at $(11,-5^7)$ and at $(-4/11,2500000/1331)$, a pole at $P_2$, and the stated value at $(11,5^7)$ (\path{u34-b3-fields.gp}). By \cref{eq:u34nu} the cusps are the points with $\nu\in\{2,\infty\}$, and $j_{\HH}(Q)\notin\{0,1,\infty\}$. Finally $\iota$ maps the points over $(11,-5^7)$ to the points over $(11,5^7)$.

(iv) \emph{The marked cover and the conic of $Q$.} Write $Z_Q$ as $v^2=r(u)$ with $r=r_2u^2+r_1u+r_0\in\Q[u]$. The branch points of class $2A$ are defined over $\Q(\sqrt{r_0})$ and those of class $3A$ over $\Q(\sqrt{r_2})$. Over $\overline{\Q}$ let $\rho_1,\rho_2$ be the roots of $r$ and let $z$ be the coordinate on $Z_Q$ with $z^2=(u-\rho_1)/(u-\rho_2)$. The four branch points are $z=\pm m$, with $m^2=\rho_1/\rho_2$, and $z=\pm1$, so that $\chi=(m-1)^2/(m+1)^2$ and
\begin{equation}\label{eq:u34conicnu}
\frac{\nu-2}{\nu+2}=\frac{4m^2}{(m^2+1)^2}=\frac{4r_0r_2}{r_1^2}.
\end{equation}
Thus the class of $r_2$ modulo squares and the value $\nu(Q)$ determine $r$ up to the substitutions $u\mapsto\lambda u$ and up to a square factor, and it remains to find the field $\Q(\sqrt{\delta})$ of the two branch points of class $3A$.

Let $\mathcal B_3\to\HH$ be the cover of \cref{lem:cuspmarked}(ii) for the class $3A$; the fiber over $Q$ is the pair of branch points of class $3A$ on $Z_Q$. Its monodromy is computed from the Hurwitz action on tuples with a marked entry (\path{u34-marked.g}): $\mathcal B_3$ is a curve of genus $11$, and $\mathcal B_3\to\HH$ is a double cover branched at the ten fixed points of $\gamma_1$ and at the two cusps of width two. Let
\[
\phi_3=-\frac{s_0\,\bigl((w^2+4)+(w+4)Z\bigr)}{\nu+2},\qquad s_0=(29w+56)\bigl(V+(35w+240)(w^2+4)\bigr),
\]
a function on $\HH$ defined over $\Q$. One has
\begin{align*}
(35w+240)^2(w^2+4)-\Psi(w)&=12\,(w-11)^2(3w-8)(29w+56),\\
\bigl((w^2+4)+(w+4)Z\bigr)\bigl((w^2+4)-(w+4)Z\bigr)&=\tfrac45\,(w^2+4)(w-1)^2\\
&=\bigl(2(w-1)Z\bigr)^2 .
\end{align*}
By the first identity $s_0$ has odd order on $\mathcal C$ exactly at the two Weierstrass points with $w^2+4=0$ and at $P_{26}$ and $P_{30}$. By the second, $(w^2+4)+(w+4)Z$ has odd order on $\HH$ exactly at the four points with $Z=0$, which are the points over these Weierstrass points. The function $\nu+2$ has five simple zeros, $17$ double zeros and the poles $P_k$. Hence $\phi_3$ has odd order exactly at the ten points of $\HH$ over the simple zeros of $\nu+2$ and at the two points over $P_2$, and $\HH(\sqrt{\phi_3})\to\mathbf P^1_\nu$ is unramified outside $-2$, $2$ and $\infty$. The script \path{u34-b3-verify.py} follows the $156$ triples $(w,Z,\sqrt{\phi_3})$ around $T=0$ and $T=1$, with certified enclosures for the roots and for the two square roots, in $3207$ and $5700$ steps. The resulting permutations have the cycle types $2^{78}$ and $1^8\,3^{16}\,5^{20}$, their product has the type $2^2\,10^4\,13^4\,15^4$, and they are simultaneously conjugate, by exactly four permutations, to the action of $(e,t)$ on the tuples of $\Sigma_0$ with a marked entry of class $3A$. So $\HH(\sqrt{\phi_3})$ and $\mathcal B_3$ are isomorphic over $\mathbf P^1_\nu$, with automorphism groups of order four. By the second identity $\phi_3\,\tau^*\phi_3=g^2$ with $g\in\Q(\HH)$ and $\tau^*g=-g$, so the lifts of $\tau$ to $\HH(\sqrt{\phi_3})$ have order four; the group is cyclic, and $\HH$ is the quotient by its only involution. An isomorphism therefore induces the identity or $\tau$ on $\HH$, and $\overline{\Q}(\mathcal B_3)$ is $\overline{\Q}(\HH)(\sqrt{\phi_3})$ or $\overline{\Q}(\HH)(\sqrt{\tau^*\phi_3})$. These fields coincide, and since $\phi_3\,\tau^*\phi_3$ is a square in $\Q(\HH)$,
\[
\Q(\mathcal B_3)=\Q(\HH)\bigl(\sqrt{c\,\phi_3}\bigr)
\]
for a squarefree integer $c$. Again the monodromy leaves a quadratic twist undetermined, and here it decides the field $\Q(\sqrt\delta)$ of the two branch points of class $3A$ at $Q$, and with it the conic. It is found from the cusps of $\mathcal B_3$ (\cref{tab:cusps}). Over each of $P_{20}$, $P_{30}$ and $P_{26}$ the curve $\mathcal B_3$ has four cusps, and their classes $\mu$ are the four classes of elements of order $10$, $15$ and $13$, respectively (\path{u34-marked.g}). A class of elements of order $10$ is fixed by no power map other than the first; a class of elements of order $15$ is fixed exactly by the first and the eleventh power, and one of order $13$ exactly by the powers $1,3,9$ (\path{u34-classfields.g}). By \cref{lem:cuspmarked}(ii) the residue fields of these cusps are $\Q(\zeta_5)$, $\Q(\zeta_5)$ and the quartic subfield of $\Q(\zeta_{13})$. On the other hand $s_0/(\nu+2)$ has the orders $10$, $16$ and $14$ at the three points of $\mathcal C$, and its leading coefficients in the parameter $w-w(P)$ are the squares $2^{-14}3^{22}5^{-8}$, $2^{-10}3^{26}5^{-20}$ and $2^{-2}3^{6}5^{-32}13^{-12}29^{24}$, while $(w^2+4)+(w+4)Z$ takes the values $20(5\pm2\sqrt5)/9$ over $w=8/3$ and $100(65\pm6\sqrt{13})/841$ over $w=-56/29$. The residue fields of $\Q(\HH)(\sqrt{c\,\phi_3})$ at the points over $P_{20}$, $P_{30}$, $P_{26}$ are therefore
\[
\Q\bigl(\sqrt5,\sqrt{-c\,(5+2\sqrt5)}\bigr),\quad\Q\bigl(\sqrt5,\sqrt{-c\,(5+2\sqrt5)}\bigr),\quad\Q\bigl(\sqrt{13},\sqrt{-c\,(65+6\sqrt{13})}\bigr).
\]
For $c=1$ these are $\Q(\zeta_5)$, $\Q(\zeta_5)$ and the quartic subfield of $\Q(\zeta_{13})$ (\path{u34-b3-fields.gp}).
\begin{table}[ht]
\centering\small
\begin{tabular}{@{}lccll@{}}
\toprule
point & order & powers & predicted field & field of the model\\
\midrule
$P_{20}$ & $10$ & $1$ & $\Q(\zeta_5)$ & $\Q\bigl(\sqrt5,\sqrt{-c\,(5+2\sqrt5)}\bigr)$\\
$P_{30}$ & $15$ & $1,\,11$ & $\Q(\zeta_5)$ & $\Q\bigl(\sqrt5,\sqrt{-c\,(5+2\sqrt5)}\bigr)$\\
$P_{26}$ & $13$ & $1,\,3,\,9$ & quartic in $\Q(\zeta_{13})$ & $\Q\bigl(\sqrt{13},\sqrt{-c\,(65+6\sqrt{13})}\bigr)$\\
\bottomrule
\end{tabular}
\caption{The four cusps of $\mathcal B_3$ over each of three points of $\mathcal C$: the order of the class $\mu$, the powers which fix it, the residue field predicted by \cref{lem:cuspmarked}(ii), and the residue field of $\Q(\HH)(\sqrt{c\,\phi_3})$. They agree only for $c=1$.}\label{tab:cusps}
\end{table}

 In general they are these fields only if $c$ is a square in $\Q(\sqrt5)$ and in $\Q(\sqrt{13})$, that is, only if $c=1$. At the two points over $(11,5^7)$ one has $s_0=2\cdot3\cdot5^{10}$, $(w^2+4)+(w+4)Z=200$ or $50$, and $\nu+2=2^2\,5^4\cdot131\cdot311^2/3^{19}$. Hence
\[
\delta\equiv\phi_3(Q)\equiv-6\cdot2\cdot393\equiv-131\pmod{\Q^{\times2}} .
\]
With $r_2=-131$ and $(\nu(Q)-2)/(\nu(Q)+2)=2^3\,61^5/(5^4\cdot131\cdot311^2)$, the relation \cref{eq:u34conicnu} gives $r_1=-5^2\cdot131\cdot311$ and $r_0=-2\cdot61^5$ for a suitable $\lambda$. The discriminant of $r$ is $3^{19}\cdot131$, and the point is verified directly.

(v) \emph{The realizations.} By (iv) the conic $Z_Q$ is isomorphic to $\mathbf P^1_{\Q}$, so \cref{prop:points}(ii) gives a regular realization with point $Q$, unique up to equivalence; its branch points of class $2A$ lie over $u=0$, hence are defined over $\Q(\sqrt{r_0})$. A regular realization of this type over $\Q(t)$ is a rational point of $\HH^{\circ}$ by \cref{prop:points}(i), since the Nielsen class is a single orbit. By (iii) this point is $Q$ or $\tau(Q)$, and the realizations with point $\tau(Q)$ are obtained from those with point $Q$ by replacing $\iota$ by $\iota\circ\varphi$, because $\tau$ is induced by $\varphi$.

(vi) This is \cref{prop:u34pullback}(ii) for the point $Q$; the branch points of $Z_Q\to\mathbf P^1_u$ are the roots of $r$.
\end{proof}

In terms of Hilbert symbols, $Z_Q$ has a rational point because $(-131,393)$ is trivial at $2$, $3$, $131$ and at infinity. The primes $61$ and $131$ in the two quadratic fields are the ones to be expected. Since $\nu(Q)-2$ is divisible by $61^5$, the point $Q$ meets modulo $61$ a cusp with $\nu=2$; there the cover of $\HH$ given by a marked point of class $2A$ is ramified and $\mathcal B_3$ is not. Since $\nu(Q)+2$ is divisible by $131$ exactly once, $Q$ meets modulo $131$ a fixed point of $\gamma_1$, where the converse holds.

By \cref{thm:u34}(v) this proves the existence of a regular realization of $U_3(4)$ of type $(2A,2A,3A,3A)$ over $\Q(t)$. The proof uses the certified continuations in (i), (ii) and (iv), exact arithmetic, and \cref{lem:cuspmarked}, which enters in (ii) and in (iv).

\begin{remark}[Further rational points of the quotient]\label{rem:u34quotient}
The curve $\mathcal C$ has at least eight rational points, those with 
\[
w\in    \left\{11,-  \frac 4 {11},\frac 8 3,- \frac {56} {29} \right\}.
\]
 Six of them are cusps and one is the image of $Q$. The eighth, 
 \[
 \left(-  \frac {56}  {29},\,   -   \frac{32500000} {24389} \right),
 \]
  is not a cusp; there $\nu=-170039803/2^{16}$, and the two points of $\HH$ over it are conjugate over $\Q(\sqrt{13})$. We have not determined $\mathcal C(\Q)$.
\end{remark}

\section{The polynomial}\label{sec:polynomial}

The extension is given by an explicit polynomial. The following theorem is proved in this section by a verification on the polynomial, in which nothing from \cref{sec:hurwitz} is used.

\begin{theorem}\label{thm:poly}
Let
\begin{align*}
f(X)={}&X^8-7328140\,X^6-1127134480\,X^5+17614872720798\,X^4\\
&+13795100796422752\,X^3-3574313729978953244\,X^2\\
&-6493268814482121098448\,X-1676431507691583902197111,
\end{align*}
and let $B$, $A$, $N$, $\tilde a\in\Z[X]$ be the polynomials of degrees $20$, $32$, $61$, $65$ recorded in the file \path{u34-poly.gp}; $B$ and $N$ are monic, $A$ has leading coefficient $3\cdot61^2$, and $\tilde a$ has leading coefficient $-1$. Put $e=131$, $c=(X-e)^2B^3$, $D=2^9\cdot131^2$,
\[
F(X,u)=D\,c\,u^2-2\,\tilde a\,u-A^2,\qquad r(u)=131u^2+1018525\,u+2\cdot61^5 .
\]
\begin{enumerate}[label=\textup{(\roman*)}]
\item One has $\tilde a^2-N^2f=-D\,c\,A^2$. The curve $Y\colon y^2=f(X)$ has genus three, and the function $u=(\tilde a+Ny)/(Dc)$ on it has degree $65$. The polynomial $F$ has degree $65$ in $X$, with leading coefficient $2u$, it is irreducible over $\Q(u)$, and $F/(2u)$ is the minimal polynomial of $X$ over $\Q(u)$. The map $u\colon Y\to\Pone$ is ramified exactly over $0$, $\infty$ and the two roots of $r$, with ramification $1\cdot2^{32}$, $1^2\,3^{21}$, $1^5\,2^{30}$ and $1^5\,2^{30}$.
\item The Galois group of $F$ over $\Q(u)$ is $G_2$, in its action on $65$ points, and the splitting field is regular over $\Q$. The inertia groups over $0$, over $\infty$ and over the roots of $r$ are generated by elements of the classes $2A$, $3A$ and $2B$. The fixed field of $U_3(4)$ is $\Q(u)\bigl(\sqrt{-r(u)}\bigr)$.
\item Let
\[
u(t)=-\frac{17371\,t^2+78732\,t+4854074}{7\,(t^2+131)} .
\]
Then $-r(u(t))$ is a square in $\Q(t)$, the polynomial $F(X,u(t))\in\Q(t)[X]$ has Galois group $U_3(4)$ over $\Q(t)$, and its splitting field is regular over $\Q$. This extension of $\Q(t)$ is ramified exactly at the two roots of $17371\,t^2+78732\,t+4854074$, which generate $\Q(\sqrt{-122})$, with inertia groups generated by involutions, and at the two roots of $t^2+131$, with inertia groups generated by elements of order three.
\end{enumerate}
\end{theorem}

The coefficients of $F$ are integers with at most $193$ digits. Part (i) is proved in \cref{sec:curve}. The Galois group of $F$ over $\overline{\Q}(u)$ is computed in \cref{sec:monodromy}, and \cref{sec:arithmetic} passes from $\overline{\Q}(u)$ to $\Q(u)$ and completes the proof.

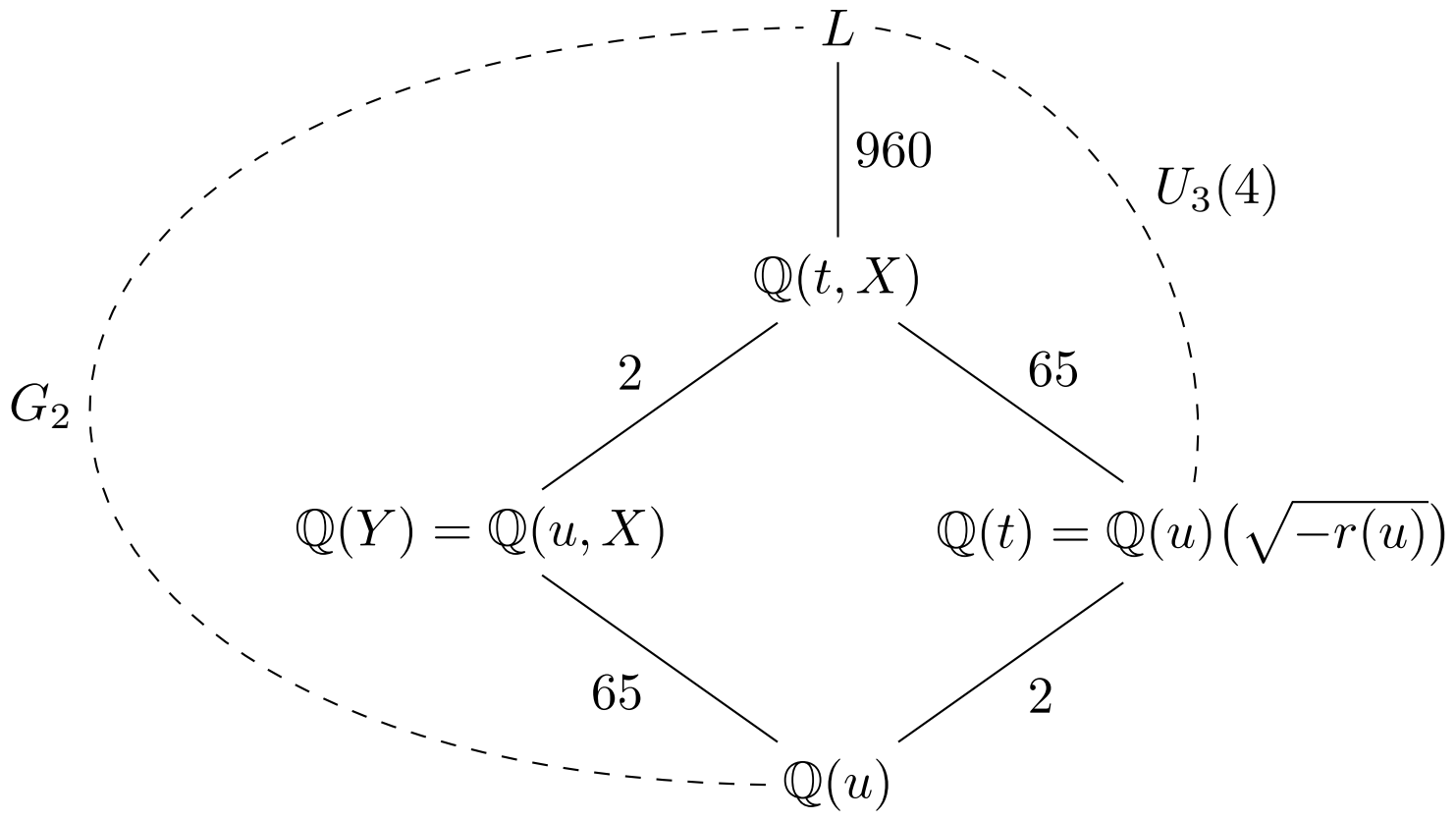
\begin{figure}[ht]
\centering
\begin{tikzpicture}[x=1cm,y=1cm,every node/.style={inner sep=3pt}]
\node (K) at (0,0) {$\Q(u)$};
\node (M) at (-2.4,1.7) {$\Q(Y)=\Q(u,X)$};
\node (T) at (2.4,1.7) {$\Q(t)=\Q(u)\bigl(\sqrt{-r(u)}\bigr)$};
\node (MT) at (0,3.4) {$\Q(t,X)$};
\node (L) at (0,5.1) {$L$};
\draw (K) -- node[below left] {$65$} (M);
\draw (K) -- node[below right] {$2$} (T);
\draw (M) -- node[above left] {$2$} (MT);
\draw (T) -- node[above right] {$65$} (MT);
\draw (MT) -- node[right] {$960$} (L);
\draw[dashed] (T.north) .. controls (2.6,3.4) and (1.6,4.9) .. node[right] {$\ U_3(4)$} (L.east);
\draw[dashed] (K.west) .. controls (-6.6,0.1) and (-6.6,5.0) .. node[left] {$G_2$} (L.west);
\end{tikzpicture}
\caption{The fields of \cref{thm:poly}. The solid lines are labelled by the degrees of the extensions, the dashed lines by their Galois groups; $L$ is the splitting field of $F$ over $\Q(u)$.}\label{fig:fields}
\end{figure}

\Cref{fig:fields} shows the fields which occur. Let $L$ be the splitting field of $F$ over $\Q(u)$. The field $\Q(Y)=\Q(u,X)$ generated by one root has degree $65$ over $\Q(u)$; it is the fixed field of the stabilizer in $G_2$ of a point of $\Omega$, a group of order $1920$. The proof goes up the other side of the diagram. \Cref{sec:monodromy,sec:arithmetic} show that $\Gal(L/\Q(u))=G_2$, that $\Q$ is algebraically closed in $L$, and that the fixed field of $U_3(4)$ is $\Q(u)(\sqrt{-r(u)})$, the function field of the conic $v^2=-r(u)$. The rational point of the conic makes this field the rational function field $\Q(t)$, and then $\Gal(L/\Q(t))=U_3(4)$ by Galois theory. Since $U_3(4)$ is transitive on $\Omega$, the field $\Q(t,X)$ has degree $65$ over $\Q(t)$, and $F(X,u(t))$ is irreducible over $\Q(t)$.

\subsection{The curve \texorpdfstring{$Y$}{Y} and the function \texorpdfstring{$u$}{u}}\label{sec:curve}

We prove \cref{thm:poly}(i). The identities and factorizations which are used are verified in exact arithmetic with PARI/GP \cite{PARI} (\path{u34-poly-verify.gp}).

The polynomial $f$ is irreducible, so $Y$ has genus three, and it has two points $\infty_\pm$ at infinity, at which $y/X^4$ takes the values $\pm1$. The norm of $\tilde a+Ny$ to $\Q(X)$ is
\[
\tilde a^2-N^2f=-D\,c\,A^2,
\]
of degree $126$. The polynomials $A$ and $B$ are squarefree and prime to each other, to $\tilde a$, to $N$ and to $X-e$, while $\tilde a(e)=N(e)=0$ and $\tilde a'(e)\ne0$. At a root of $B$ the norm vanishes to order three and $\tilde a$ does not vanish, so $\tilde a+Ny$ vanishes to order three at one of the two points over it and $u$ has a triple pole at the other. At the two points over $X=e$ the function $u$ has simple poles. Since $\tilde a$ and $N$ have the leading coefficients $-1$ and $1$, the function $\tilde a+Ny$ has a pole of order $65$ at $\infty_-$, where $u$ has a pole of order three; comparing degrees, $u$ has a simple zero at $\infty_+$. Thus $u$ has $2+60+3=65$ poles, counted with multiplicity. In the same way its zeros are $\infty_+$ and one double zero over each root of $A$. This gives the ramification $1\cdot2^{32}$ over $0$ and $1^2\,3^{21}$ over $\infty$.

By the norm identity $u\,(Dc\,u-2\tilde a)=(\tilde a+Ny)(Ny-\tilde a)/(Dc)=A^2$, that is, $F(X,u)=0$ on $Y$. Moreover $y=(Dc\,u-\tilde a)/N$, so $u$ and $X$ generate the function field of $Y$, which has degree $65$ over $\Q(u)$. The coefficient of $X^{65}$ in $F$ is that of $-2\tilde a\,u$, which is $2u$. Hence $F$ has degree $65$ in $X$, it is irreducible over $\Q(u)$, and $F/(2u)$ is the minimal polynomial of $X$ over $\Q(u)$. The function $X$ has its poles at $\infty_\pm$, which lie over $u=0$ and $u=\infty$. Therefore, for every $u_0\ne0$ in $\overline{\Q}$,
\begin{equation}\label{eq:fibre}
F(X,u_0)=2u_0\prod_{P}\bigl(X-X(P)\bigr)^{e_P},
\end{equation}
where $P$ runs over the points of $Y$ with $u(P)=u_0$ and $e_P$ is the ramification index of $u$ at $P$.

Let $\rho_1$, $\rho_2$ be the roots of $r$. The polynomial $r$ is irreducible over $\Q$, so $F(X,\rho_1)$ and $F(X,\rho_2)$ are conjugate over $\Q$. They have no common root. Indeed, let $\alpha u+\beta$, with $\alpha,\beta\in\Z[X]$, be the remainder of $131\,F$ in the division by $r(u)$ as polynomials in $u$. If $F(x_0,\rho_1)=F(x_0,\rho_2)=0$, then $F(x_0,u)$, a polynomial of degree at most two in $u$, is divisible by $r(u)$, so $\alpha(x_0)=\beta(x_0)=0$; but $\alpha$ and $\beta$ are coprime. The product $F(X,\rho_1)\,F(X,\rho_2)$ is, up to a constant, $S\,T^2$ with $S,T\in\Q[X]$ squarefree and coprime, of degrees $10$ and $60$, and $T$ is prime to $N$. Since the two factors have no common root and are conjugate, $F(X,\rho_i)$ is, up to a constant, $S_i\,T_i^2$ with $S=S_1S_2$, $T=T_1T_2$, and with $S_i$ of degree $5$ and $T_i$ of degree $30$. By \eqref{eq:fibre} a root $x_0$ of $T_i$ is the value of $X$ at one point of the fiber over $\rho_i$ with ramification index two, or at two unramified points of this fiber. Two points of $Y$ with the same value $x_0$ of $X$ are $(x_0,\pm y_0)$ with $y_0\ne0$, and they have the same value of $u$ only if $N(x_0)=0$, which is excluded. Hence $u$ has $30$ points of index two and five unramified points over each root of $r$. The ramification found has total index
\[
32+42+30+30=134=2\cdot65+2\cdot3-2,
\]
so by the Riemann--Hurwitz formula there is no other.

\begin{remark}
The hyperelliptic involution of $Y$ exchanges the two simple poles of $u$ and maps its simple zero to a triple pole. The coordinate $X$ is normalized as follows: $X=3^2\,61^4\,x/2$, where $x$ is the function of degree two on $Y$ with $x=1/u+O(1)$ at the simple zero of $u$ and for which $f$ has no term of degree seven. The discriminant of $f$ is divisible by no primes below $10^6$ other than $2$, $3$, $61$ and $131$.
\end{remark}

\subsection{The monodromy}\label{sec:monodromy}

Let $\rho_1\approx-2398.2$ and $\rho_2\approx-5376.8$ be the roots of $r$. By \cref{thm:poly}(i) the branch points of $u$ are $0$, $\rho_1$, $\rho_2$ and $\infty$. We compute the action of the fundamental group of $\C\setminus\{0,\rho_1,\rho_2\}$ on the $65$ roots of $F(X,u)$.

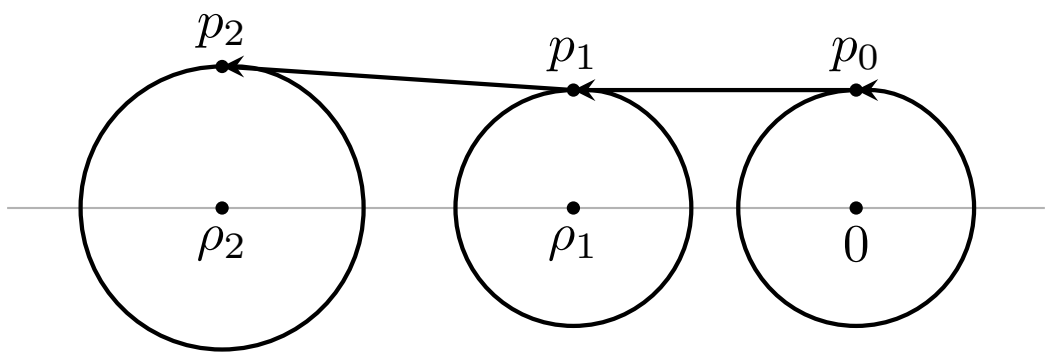
\begin{figure}[ht]
\centering
\begin{tikzpicture}[scale=0.78,>=stealth]
\draw[gray!60] (-7.2,0)--(1.6,0);
\foreach \x/\r/\lab in {-5.377/1.2/{\rho_2},-2.398/1.0/{\rho_1},0/1.0/{0}} {
  \draw[thick,->] (\x,\r) arc (90:450:\r);
  \fill (\x,0) circle (1.6pt) node[below] {$\lab$};
}
\draw[thick] (-2.398,1.0)--(0,1.0);
\draw[thick] (-2.398,1.0)--(-5.377,1.2);
\fill (-2.398,1.0) circle (1.6pt) node[above] {$p_1$};
\fill (0,1.0) circle (1.6pt) node[above] {$p_0$};
\fill (-5.377,1.2) circle (1.6pt) node[above] {$p_2$};
\end{tikzpicture}
\caption{The loops in the $u$-plane, with base point $p_1$. One unit is $1000$.}\label{fig:loops}
\end{figure}

Let $p_1=\rho_1+1000\,\sqrt{-1}$, $p_0=1000\,\sqrt{-1}$ and $p_2=\rho_2+1200\,\sqrt{-1}$. The base point is $p_1$. The first loop is the circle about $\rho_1$ through $p_1$. The second is the segment from $p_1$ to $p_0$, the circle about $0$ through $p_0$, and the segment back; the third is the segment from $p_1$ to $p_2$, the circle about $\rho_2$ through $p_2$, and the segment back. All circles are run counterclockwise. The three closed discs are pairwise disjoint and each contains exactly one of the points $\rho_1$, $0$, $\rho_2$. The segment from $p_1$ to $p_0$ is tangent to two of the circles, and the segment from $p_1$ to $p_2$ meets the disc about $\rho_2$ in a chord which cuts off a part not containing $\rho_2$. The classes of the three loops therefore generate the fundamental group of $\C\setminus\{0,\rho_1,\rho_2\}$ with base point $p_1$, which maps onto that of the complement of the four branch points in $\Pone$ (\cref{fig:loops}).

\begin{proposition}\label{prop:monodromy}
Let $g_1$, $g_0$, $g_2$ be the permutations of the $65$ roots of $F(X,p_1)$ along the loops about $\rho_1$, $0$ and $\rho_2$.
\begin{enumerate}[label=\textup{(\roman*)}]
\item Their cycle types are $1^5\,2^{30}$, $1\cdot2^{32}$ and $1^5\,2^{30}$, and the product of the three in any order has the type $1^2\,3^{21}$.
\item They generate a group $H$ of order $124800$, and there is a bijection from the roots to $\Omega$ which carries $H$ to $G_2$. The commutator subgroup of $H$ has order $62400$; it contains $g_0$ and does not contain $g_1$ or $g_2$.
\item The Galois group of $F$ over $\overline{\Q}(u)$ is $H$, and the normalizer of $H$ in the symmetric group of the roots is a group of order $249600$ which corresponds to $\Aut(U_3(4))$.
\end{enumerate}
\end{proposition}

\begin{proof}
The roots of $F$ were followed by the certified procedure of \cref{lem:rouche}, in the variable $X/500$ and with $600$ bits (\path{u34-poly-monodromy.py}). The continuation was certified along the three circles in $2510$, $3893$ and $1142$ steps and along the two segments in $1168$ and $663$ steps, $9376$ steps in all, and the permutations are recorded in \path{u34-poly-monodromy.g}. Statements (i) and (ii) are computed from them (\path{u34-poly-group.g}); the bijection is found by GAP as a conjugating element in $S_{65}$. The Galois group of $F$ over $\C(u)$, which is that over $\overline{\Q}(u)$, is the image of the fundamental group, and the normalizer is given by \cref{lem:group}(b).
\end{proof}

By (ii) and \cref{lem:group}(a), the inertia groups over $0$, $\infty$ and the roots of $r$ are generated by elements of the classes $2A$, $3A$ and $2B$, in accordance with the ramification found in \cref{sec:curve}.

\subsection{The field of constants and the quadratic subfield}\label{sec:arithmetic}

The passage from $\overline{\Q}(u)$ to $\Q(u)$ uses specializations and their reductions modulo primes, through the following standard facts.

\begin{lemma}\label{lem:constants}
Let $F\in\Z[u][X]$ have degree $n$ in $X$ and be irreducible over $\overline{\Q}(u)$. Let $\Gamma\le S_n$ be its Galois group over $\Q(u)$, $H$ its Galois group over $\overline{\Q}(u)$, and $k$ the field of constants of its splitting field. Then $H$ is the Galois group of $F$ over $k(u)$, it is normal in $\Gamma$, and $\Gamma/H\simeq\Gal(k/\Q)$. Let $u_0\in\Q$ be such that $F(X,u_0)$ is separable of degree $n$, and let $F_{u_0}\in\Z[X]$ be the primitive polynomial proportional to it.
\begin{enumerate}[label=\textup{(\alph*)}]
\item The splitting field $E$ of $F_{u_0}$ over $\Q$ contains $k$, and its Galois group is, as a permutation group of the roots, conjugate to a subgroup of $\Gamma$, compatibly with the restriction to $k$.
\item A prime which ramifies in $E$ divides the leading coefficient or the discriminant of $F_{u_0}$. For a prime $p$ which divides neither, the degrees of the irreducible factors of $F_{u_0}$ modulo $p$ are the cycle type of an element of $\Gamma$ whose restriction to $k$ is the Frobenius automorphism at $p$.
\end{enumerate}
\end{lemma}

\begin{proof}
The first statements are Galois theory. Since $F(X,u_0)$ is separable of degree $n$, the place $u=u_0$ of $\Q(u)$ is unramified in the splitting field of $F$, the residue field at a place above it is $E$, and the decomposition group is mapped isomorphically to $\Gal(E/\Q)$, compatibly with the actions on the roots and with the restriction to $k$. This is (a). Part (b) is Dedekind's theorem on the reduction of a polynomial modulo a prime, applied to $F_{u_0}$.
\end{proof}

The factorizations which follow are therefore not statistical evidence: they decide between finitely many possibilities for $k$ and for the quadratic subfield, and the list of these possibilities is established first.

Let now $F$ be the polynomial of \cref{thm:poly}, $\Gamma$ its Galois group over $\Q(u)$ and $k$ the field of constants of its splitting field. By \cref{prop:monodromy} the group $H$ of \cref{lem:constants} is $G_2$.

\begin{proposition}\label{prop:arithmetic}
One has $\Gamma=G_2$ and $k=\Q$, and the fixed field of $U_3(4)$ in the splitting field of $F$ is $\Q(u)\bigl(\sqrt{-r(u)}\bigr)$.
\end{proposition}

\begin{proof}
The group $\Gamma$ normalizes $G_2$, so by \cref{prop:monodromy}(iii) it is $G_2$ or $\Aut(U_3(4))$, and $k$ has degree one or two over $\Q$. If $k=\Q$, the fixed field of $U_3(4)$ is a quadratic extension of $\Q(u)$ which is ramified exactly at the roots of $r$, because the inertia groups there are generated by elements outside $U_3(4)$ and those over $0$ and $\infty$ lie in $U_3(4)$. It is therefore $\Q(u)(\sqrt{\kappa\,r(u)})$ for a squarefree integer $\kappa$.

For an integer $u_0$ such that $F(X,u_0)$ is separable of degree $65$, let $F_{u_0}$ be as in \cref{lem:constants}, with leading coefficient $\ell(u_0)$, and let $E$ be its splitting field. By \cref{lem:constants}(a) the field $E$ contains $k$, and it contains $\sqrt{\kappa\,r(u_0)}$ when $k=\Q$.

A prime which ramifies in $k$, or which divides $\kappa$, ramifies in $E$ for every such $u_0$ with $r(u_0)$ prime to it, and so it divides $\ell(u_0)$ or the discriminant of $F_{u_0}$ by \cref{lem:constants}(b). The greatest common divisor of the numbers
\[
u_0\,r(u_0)\,\ell(u_0)\,\operatorname{disc}F_{u_0},\qquad u_0=1,\ldots,7,
\]
is divisible only by the primes $2$, $3$ and $131$ (\path{u34-poly-verify.gp}). So $k$ is $\Q$ or one of the $15$ quadratic fields unramified outside these primes, and $\kappa$ is one of the $16$ squarefree integers divisible only by them.

For the $54$ primes $p$ between $1009$ and $1399$ and for $1\le u_0<60$, the factorizations of $F(X,u_0)$ modulo $p$ were computed (\path{u34-poly-check.py}). The $3179$ separable ones all have the cycle type of an element of $G_2$, and such a type occurs for each of the $54$ primes. If $k$ were quadratic, then by \cref{lem:constants}(b) these factorizations would give, at a prime inert in $k$, elements of $\Gamma$ which act nontrivially on $k$ and so lie outside $G_2$; by \cref{lem:group}(a) their cycle types would not be types of $G_2$. Hence all $54$ primes would split in $k$. But none of the $15$ fields is split at all $54$ primes. Therefore $k=\Q$ and $\Gamma=G_2$.

Finally, let $k=\Q$. Since $F(X,\rho_i)$ has multiple roots, the discriminant of $F$ with respect to $X$ vanishes at the roots of $r$; so it is divisible by $r(u)$ in $\Q[u]$, and in $\Z[u]$ because $r$ is primitive. Hence, if $F(X,u_0)$ is separable of degree $65$ modulo one of the $54$ primes $p$, then $p$ does not divide $r(u_0)$, and $p$ is unramified in $\Q(\sqrt{\kappa\,r(u_0)})$. Again by \cref{lem:group}(a), the cycle type of the factorization modulo $p$ is then that of an element of $U_3(4)$ exactly when the Frobenius automorphism at $p$ fixes $\sqrt{\kappa\,r(u_0)}$, that is, when $\kappa\,r(u_0)$ is a square modulo $p$. Among the $16$ values of $\kappa$ only $\kappa=-1$ is compatible with the $3179$ factorizations; each of the others disagrees with at least $1355$ of them.
\end{proof}

\begin{proof}[Proof of \cref{thm:poly}]
Part (i) is proved in \cref{sec:curve}, and (ii) is \cref{prop:arithmetic} together with \cref{prop:monodromy}. For (iii), one has $-r(u(t))=v(t)^2$ with
\[
v(t)=\frac{39366}{7}+t\,\Bigl(u(t)+\frac{17371}{7}\Bigr);
\]
this is the parametrization of the conic $v^2=-r(u)$ by the lines through its rational point $(u,v)=(-17371/7,\,39366/7)$. Since $u(t)$ has degree two, $\Q(t)=\Q(u)(\sqrt{-r(u)})$, and by (ii) the splitting field of $F(X,u(t))$ over $\Q(t)$ is the splitting field of $F$ over $\Q(u)$, with group $U_3(4)$ over $\Q(t)$; it is regular over $\Q$. The map $t\mapsto u(t)$ is ramified exactly over the roots of $r$, where the inertia groups of the splitting field of $F$ have order two and are not contained in $U_3(4)$; so the extension of $\Q(t)$ is unramified there. It is ramified at the two points over $u=0$ and at the two points over $u=\infty$, with the inertia groups of (ii). These are the roots of the numerator and of the denominator of $u(t)$. The discriminant of $17371\,t^2+78732\,t+4854074$ is $-2^3\cdot7^2\cdot61^5$.
\end{proof}

\section{The main theorem}\label{sec:main}

\begin{theorem}\label{thm:main}
There is a Galois extension of $\Q(t)$ with group $U_3(4)$, regular over $\Q$, which is ramified at exactly four points: the two roots of $17371\,t^2+78732\,t+4854074$, with inertia groups of order two, and the two roots of $t^2+131$, with inertia groups of order three. Every regular Galois extension of $\Q(t)$ with group $U_3(4)$ which is ramified at exactly four points, with inertia groups of orders $2$, $2$, $3$, $3$, is obtained from it by an automorphism of $\Q(t)$ over $\Q$. In particular $U_3(4)$ is the Galois group of infinitely many Galois extensions of $\Q$.
\end{theorem}

\begin{proof}
The splitting field of $F(X,u(t))$ over $\Q(t)$ has the properties of the first statement by \cref{thm:poly}(iii), and the last statement follows from the first by Hilbert's irreducibility theorem. Let $F'$ be a regular Galois extension of $\Q(t)$ with group $U_3(4)$ which is ramified at exactly four points with inertia groups of orders $2$, $2$, $3$, $3$. Its inertia classes are $2A$, $2A$, $3A$, $3A$, since these are the only classes of elements of order two and three. With an isomorphism $\iota$ of $U_3(4)$ with its Galois group it is a regular realization of type $\mathbf C$, and by \cref{prop:points}(i) it determines a point of $\HH^{\circ}(\Q)$. By \cref{thm:u34}(iii) this point is $Q$ or $\tau(Q)$, and since $\tau$ is induced by $\varphi$, we may replace $\iota$ by $\iota\circ\varphi$ and assume that it is $Q$. The same applies to the splitting field $L$ of $F(X,u(t))$: with a suitable isomorphism of $U_3(4)$ with its Galois group it is a regular realization of type $\mathbf C$ with the point $Q$. The point $Q$ lies in $\HH^{\mathrm f}$, and by \cref{prop:points}(ii) the regular realizations of type $\mathbf C$ over $\Q$ with the point $Q$ form one equivalence class. Hence $F'$ is obtained from $L$ by an element of $\operatorname{PGL}_2(\Q)$ acting on $t$.
\end{proof}

In particular the extension of \cref{thm:poly}(iii) is the one of \cref{thm:u34}(v). The existence statement of \cref{thm:main} is thus proved twice, by \cref{thm:u34} and by \cref{thm:poly}, and \cref{prop:arithmetic} confirms the conic of \cref{thm:u34}(iv) by an argument which does not use \cref{lem:cuspmarked}.

A polynomial over $\Q$ with the group $U_3(4)$ can be written down: it is the specialization of $F(X,u(t))$ at the point $t=\infty$ of the $t$-line, where $u=-17371/7$.

\begin{corollary}\label{cor:Q}
The polynomial
\[
P(X)=7^2\,F\bigl(X,-\tfrac{17371}{7}\bigr)=17371^2\,D\,c+2\cdot7\cdot17371\,\tilde a-7^2A^2\ \in\Z[X]
\]
is irreducible of degree $65$, and its Galois group over $\Q$ is $U_3(4)$.
\end{corollary}

\begin{proof}
Put $u_0=-17371/7$. The polynomial $P$ has degree $65$ and the leading coefficient $-2\cdot7\cdot17371$, and its coefficients have no common factor (\path{u34-poly-Q.gp}). Modulo $5$ it is the product of $13$ distinct irreducible factors of degree $5$, and modulo $17$ it is the product of five distinct irreducible factors of degree $13$; neither prime divides the leading coefficient. Hence the degree of every factor of $P$ over $\Q$ is divisible by $5$ and by $13$, so $P$ is irreducible, and it is separable. By \cref{lem:constants}(a) its Galois group is a decomposition group of the splitting field of $F$ at a place over $u=u_0$. Since $-r(u_0)=(39366/7)^2$ is a nonzero square, this place splits in $\Q(u)(\sqrt{-r(u)})$, so by \cref{thm:poly}(ii) the Galois group of $P$ is a subgroup of $U_3(4)$. By \cref{lem:constants}(b) the two factorizations give elements of order $5$ and of order $13$ in it, and it is $U_3(4)$ by \cref{lem:group}(c).
\end{proof}

The polynomial $P$ is not monic, and its coefficients have up to $194$ digits. \Cref{sec:appendix} gives a monic polynomial $P_0$ of degree $65$, with coefficients of at most $127$ digits, which is obtained from $P$ by a change of the generator of the field $\Q[X]/(P)$ and therefore has the same Galois group.

\section{Three branch points}\label{sec:rigidity}

A realization with three branch points would be simpler, but the rigidity method does not provide one for $U_3(4)$. This section makes the statement precise; it is not used elsewhere. For a group $G$ with trivial center, consider the triples of conjugacy classes whose class vector is stable under all power maps up to a permutation of its entries; these are the class vectors to which the rationality criteria apply. Such a vector is \emph{rigid} if the triples $(x,y,z)$ in the three classes with $xyz=1$ which generate $G$ form one orbit under conjugation.

\begin{proposition}\label{prop:rigid}
Consider the class vectors of three classes which are stable under the power maps up to permutation.
\begin{enumerate}[label=\textup{(\roman*)}]
\item In $U_3(4)$ and in $G_2$ none of them is rigid: a vector which has generating triples has at least four, respectively two, orbits of them.
\item The group $\Aut(U_3(4))$ has four classes of elements of order $16$. Each of them is fixed by the powers congruent to $1$ modulo $4$, and the third power exchanges them in two pairs, which we denote by $16A$, $16B$ and $16C$, $16D$.
\item In $\Aut(U_3(4))$ exactly two of them are rigid, $(2A,16A,16B)$ and $(2A,16C,16D)$.
\end{enumerate}
\end{proposition}

\begin{proof}
Let $H$ be one of the three groups, let $C_1$, $C_2$, $C_3$ be conjugacy classes of $H$, and fix $x\in C_1$ and $y_0\in C_2$. Every triple in $C_1\times C_2\times C_3$ with product one is conjugate to a triple $(x,y,(xy)^{-1})$, and two such triples are conjugate exactly when the elements $y$ are conjugate under the centralizer $C_H(x)$. So the orbits of $H$ on these triples correspond to the orbits of $C_H(x)$ on the set of all $y\in C_2$ with $(xy)^{-1}\in C_3$, and the orbits of $C_H(x)$ on $C_2$ correspond to the double cosets $C_H(y_0)\,h\,C_H(x)$, by $h\mapsto h^{-1}y_0h$. The number of orbits which consist of generating triples does not depend on the order of the three classes, since $(x,y,z)\mapsto(y,y^{-1}xy,z)$ and $(x,y,z)\mapsto(x,z,z^{-1}yz)$ are bijections between the sets of triples with product one which commute with conjugation and preserve the group generated.

The program \path{u34-rigidity.g} computes with GAP \cite{GAP} the classes of $H$ and their power maps. For each pair of nontrivial classes $C_1$, $C_2$ it takes representatives $y$ of all orbits of $C_H(x)$ on $C_2$ from the double cosets, determines the class of $(xy)^{-1}$, and decides whether $x$ and $y$ generate $H$. This gives, for every multiset of three nontrivial classes, the number of orbits of generating triples. The class vectors which are stable under the power maps up to permutation and have a triple with product one number $11$ in $U_3(4)$, $34$ in $G_2$ and $196$ in $\Aut(U_3(4))$. Of these, $5$, $9$ and $30$ have generating triples, and the smallest number of orbits of generating triples is $4$, $2$ and $1$. This is (i). In $\Aut(U_3(4))$ exactly two of the $30$ vectors have one orbit. Each of them consists of the class $2A$ of the involutions of $U_3(4)$ and of two classes of elements of order $16$, and these two classes are exchanged by the third power. The same program computes all powers of the elements of order $16$, which gives (ii), and with it (iii).
\end{proof}

\begin{proposition}\label{prop:rigidreal}
Let $\mathbf v$ be one of the two rigid class vectors of $A=\Aut(U_3(4))$ of \cref{prop:rigid}.
\begin{enumerate}[label=\textup{(\roman*)}]
\item There is a Galois extension $N$ of $\Q(t)$ with group $A$, regular over $\Q$, which is ramified exactly at $t=\infty$ and $t=\pm\sqrt{-1}$, with inertia classes $\mathbf v$.
\item The fixed field of $U_3(4)$ in $N$ is the function field of a conic $\mathcal Z_0\colon x^2+y^2=\beta$ with $\beta\in\Q^\times$, so $U_3(4)$ is the Galois group of a regular extension of $\Q(\mathcal Z_0)$. The fixed field of $G_2$ is a rational function field $\Q(s)$, so $G_2$ has a regular realization over $\Q(s)$.
\item One has $\beta<0$. The conic $\mathcal Z_0$ has no real point; in particular the fixed field of $U_3(4)$ in $N$ is not a rational function field.
\end{enumerate}
\end{proposition}

\begin{proof}
(i) Let $\mathbf v=(2A,16A,16B)$; the other case is the same. By Riemann's existence theorem and rigidity there is a Galois cover $X\to\mathbf P^1$ over $\overline{\Q}$ with group $A$, ramified at $\infty$, $\sqrt{-1}$, $-\sqrt{-1}$ with the inertia classes $2A$, $16A$, $16B$, and since $Z(A)=1$ it is unique up to a unique isomorphism which commutes with $A$. Let $\sigma\in\Gal_{\Q}$. By the branch cycle argument \cite[Lemma~2.8]{Voelklein}, the conjugate cover is ramified at $\infty$ and $\pm\sqrt{-1}$, and its inertia class at $\sigma(z)$ is the class at $z$ raised to the power $n'$, where $n=\kappa_{16}(\sigma)$ and $nn'\equiv1$ modulo $16$. If $\sigma$ fixes $\sqrt{-1}$, then $n'\equiv1$ modulo $4$ and the classes are unchanged. Otherwise $n'\equiv3$ modulo $4$, the two points $\pm\sqrt{-1}$ are exchanged and so are the classes $16A$ and $16B$. In both cases the conjugate cover has the same branch points with the same classes, so it is isomorphic to $X$. Hence the cover with its action of $A$ is a point of the inner Hurwitz space of $(A,\mathbf v)$ which is fixed by $\Gal_{\Q}$, and since $Z(A)=1$ it is defined over $\Q$ together with the action of $A$ \cite[\S2.1, Corollary~1]{FV}.

(ii) Let $L$ be the fixed field of $U_3(4)$ in $N$. It is a cyclic extension of $\Q(t)$ of degree four, regular over $\Q$. The class $2A$ lies in $U_3(4)$, and an element of order $16$ generates $A$ modulo $U_3(4)$; so $L$ is totally ramified at $\pm\sqrt{-1}$ and unramified elsewhere, and it has genus zero. By Kummer theory $L(\sqrt{-1})$ is generated over $\Q(\sqrt{-1})(t)$ by a function $\theta$ with
\[
\theta^4=c_0\,\frac{t-\sqrt{-1}}{t+\sqrt{-1}},\qquad c_0\in\Q(\sqrt{-1})^\times .
\]
Let $\rho$ be the automorphism of $L(\sqrt{-1})$ over $L$ of order two. Then $\rho(\theta)^4=\rho(c_0)c_0/\theta^4$, so $\rho(\theta)=\beta/\theta$ with $\beta\in\Q(\sqrt{-1})$, and $\rho^2=1$ gives $\rho(\beta)=\beta$, that is, $\beta\in\Q^\times$. The functions $x=(\theta+\beta/\theta)/2$ and $y=(\theta-\beta/\theta)/(2\sqrt{-1})$ are fixed by $\rho$, they satisfy $x^2+y^2=\beta$, and $\theta=x+y\sqrt{-1}$. Hence $L=\Q(x,y)$ is the function field of $\mathcal Z_0$. The fixed field of $G_2$ is the subfield of $L$ fixed by $\theta\mapsto-\theta$, which is generated by $x^2-y^2$ and $2xy$; it is the function field of the conic $X^2+Y^2=\beta^2$, which has the rational point $(\beta,0)$.

(iii) Let $t_0$ be a real number. Complex conjugation acts on the fiber of $X$ over $t_0$ and commutes with the action of $A$, so for a point $x_0$ of this fiber it maps $a\cdot x_0$ to $ac\cdot x_0$ for an element $c\in A$ with $c^2=1$. On the fiber of $\mathcal Z_0=X/U_3(4)$, which is a principal homogeneous space under $A/U_3(4)$, it acts by the translation by the image of $c$. Hence $\mathcal Z_0$ has a real point over $t_0$ if and only if $c\in U_3(4)$. The element $c$ is determined by the branch cycles. Let $\gamma_+$ and $\gamma_-$ be the loops based at $t_0$ which go along a path to a point near $\sqrt{-1}$, respectively along the complex conjugate path to a point near $-\sqrt{-1}$, around that point counterclockwise, and back. Complex conjugation maps $\gamma_+$ to $\gamma_-^{-1}$. If $g_+$, $g_-$ are the corresponding branch cycles in the action on $x_0$, this gives
\[
c\,g_+\,c^{-1}=g_-^{-1},\qquad c\,g_-\,c^{-1}=g_+^{-1}
\]
\cite[\S2.3, Propositions~2.2 and~2.3]{DebesFried1994}, and these equations determine $c$, because $g_+$ and $g_-$ generate $A$ and $Z(A)=1$. For both rigid class vectors the element $c$ computed from a generating triple is an involution which lies in $G_2$ and not in $U_3(4)$ (\path{u34-rigidity.g}); by rigidity this does not depend on the triple. Every real point of $\mathcal Z_0$ lies over a real point of the $t$-line, and the branch points are not real. So $\mathcal Z_0$ has no real point, and $\beta<0$.
\end{proof}

Thus rigidity gives regular realizations of $\Aut(U_3(4))$ and of $G_2$ over $\Q(t)$, and leads for $U_3(4)$ to the question of a rational point on a conic, as in the introduction; but the conic has no real point. In \cref{thm:poly} the conic is $v^2=-r(u)$. Both roots of $r$ are real and negative, and $-r$ is positive between them, so the obstruction at the real place is absent.

\section{How the equations were found}\label{sec:found}

The equations of \cref{thm:u34,thm:poly} were found numerically and then proved as above. This section describes the two computations. They are not certified, no statement of the paper depends on them, and their programs are in the directories \path{u34-compute} and \path{u34-poly-compute}.

\subsection{The Hurwitz curve}
The weight-two cusp forms on the subgroup of the modular group which belongs to $\HH$ were computed to $36$ digits, and those of the quotient by $\tau$ to $56$ digits (\path{u34-compute}). The three forms on $\HH$ satisfy one quadratic relation, so the curve is hyperelliptic, and the equation \cref{eq:u34octic} was recognized from the quotient $x$ of two forms of an echelon basis at a cusp of width $20$. The discriminant of its right-hand side is $2^{16}3^{4}5^{17}13$, and the primes dividing it are those dividing $|U_3(4)|=2^6\cdot3\cdot5^2\cdot13$. The cusps of $\mathcal C$ were recognized from the numerical data. Besides the four poles of $\nu$ they are the points $(11,-5^7)$ and $(-4/11,\,2500000/1331)$, of width one; four points of width three, namely the two Weierstrass points with $181w^2-732w+1776=0$ and the two points with $137w^2-264w+1952=0$ and $V=-(12045000w+64380000)/18769$; and five points of width five, with
\[
1807w^5+963240w^4-75280w^3+4902080w^2+706560w-1964032=0 .
\]
The conditions on the divisor of $\nu-2$ are $58$ linear equations for the $57$ coefficients of $A_\nu$ and $B_\nu$, and over $\Q$ they have a solution, unique up to a factor (\path{u34-compute/nu.gp}); the factor is fixed by the values of $\nu$ computed from the modular function $\lambda$. The function $s_0$ was found from the divisor of $\phi_3$: it has odd order exactly at the two Weierstrass points over $w^2+4=0$ and at $P_{26}$ and $P_{30}$, because $(35w+240)^2(w^2+4)-\Psi(w)=12\,(w-11)^2(3w-8)(29w+56)$.

For the two other curves of the table the same numerical computation, which has not been certified, shows that the curve of genus three with classes $(2A,2A,5E,5F)$ is a plane quartic and that the curve of genus five is not hyperelliptic.

\subsection{The polynomial}\label{sec:polyfound}
The point $Q$ gives the branch points $0$, $\infty$, $\rho_1$, $\rho_2$ of a cover of the $u$-line with group $G_2$, but not its branch cycles. The covers of degree $65$ with group $G_2$, these four branch points and the classes $2A$, $3A$, $2B$, $2B$ are described by $39$ classes of tuples of permutations, in accordance with the degree of $\nu$ on $\mathcal C$. Five of them are invariant under complex conjugation, as are five points of $\mathcal C$ over $\nu(Q)$, and following the real points of $\mathcal C$ to a cusp leaves two candidates for the cover which belongs to $Q$. By the Riemann--Hurwitz formula such a cover has genus three.

No equation was available, so the cover was computed from its branch cycles through its holomorphic differentials, in the following steps:
\[
\begin{array}{c}
\text{branch cycles}\ \to\ \text{local expansions of differentials}\ \to\ \text{numerical kernel}\\[2pt]
\to\ \text{hyperelliptic model}\ \to\ \text{refinement}\ \to\ \text{exact verification}.
\end{array}
\] A holomorphic differential on the cover is given, near each point of the $u$-line, by $65$ power series in a local parameter, one for each sheet, or for each cycle of the local monodromy at a branch point. We used six centres, the four branch points and two further points, in a coordinate on the $u$-line in which the branch points lie on the unit circle and the discs of convergence overlap. The conditions that the series at neighbouring centres agree at sample points of the overlaps, with the identification of the sheets given by the tuple, form a linear system; with series truncated at order $32$ it has $12480$ unknowns and $24960$ equations. Its numerical null space has dimension three, the genus: the three smallest singular values are of the order $10^{-10}$ and the next is $4\cdot10^{-3}$. The three differentials satisfy a quadratic relation. So the canonical image is a conic, the cover is a hyperelliptic curve $Y$, and the quotient $x$ of two differentials is a function of degree two on it.

The divisor of $u$ on $Y$ is known from the tuple, which gives the shape $y^2=f(X)$, $u=(a+by)/c$ of \cref{thm:poly} with $c$ and the norm of $a+by$ determined by the poles and the zeros of $u$. A first approximation of all polynomials was obtained from the numerical cover; the roots of $b$, which are the values of $X$ at pairs of points with the same value of $u$, were found by a search on it. The model was then refined by a Gauss--Newton iteration in which the unknowns are the roots of the polynomials and the identities are imposed at points; with the coefficients as unknowns the iteration did not converge, the roots being spread over two orders of magnitude. It converges quadratically once the coincidence of the values of $X$ at the two simple poles of $u$ is built into the unknowns, and it was carried to $1400$ digits. For one of the two candidates the coefficients are then recognized as rational numbers by continued fractions, with denominators divisible only by $2$, $3$, $61$ and $131$; for the other they are not recognized as rational numbers at $120$ digits. The identities of \cref{thm:poly}(i) were then verified exactly.

\appendix
\section{Computations}\label{sec:computations}

\subsection{The programs}\label{sec:files}

\begingroup\emergencystretch=3em
The proofs use PARI/GP~2.15.4 \cite{PARI}, GAP~4.12.1 \cite{GAP} and Python~3.13 with \texttt{python-flint}~0.9.0 (FLINT~3.6.0), whose ball arithmetic is that of \cite{Johansson2017}. \Cref{tab:files} lists, for each statement, the method and the programs; all files accompany the paper, each program has its output in a file with the extension \texttt{.log}, and the file \path{README.txt} gives the commands. The polynomials of \cref{thm:poly} are in \path{u34-poly.gp}, which defines $e$, $D$, $f$, $B$, $A$, $N$, $\tilde a$, $c$, $F$, $r$, $u(t)$ and $v(t)$ for PARI/GP, and in \path{u34-poly-data.json}; the polynomial $P$ is in \path{u34-poly-Q.txt}, and $P_0$ in \path{u34-poly-red.txt}.
\endgroup

\Needspace{0.5\textheight}
\begin{center}\footnotesize
\begin{tabular}{@{}>{\raggedright\arraybackslash}p{0.27\textwidth}>{\raggedright\arraybackslash}p{0.36\textwidth}>{\raggedright\arraybackslash}p{0.30\textwidth}@{}}
\toprule
statement & method & programs\\
\midrule
\Cref{lem:group} & finite groups (GAP) & \path{u34-cycletypes.g}, \path{u34-maxes.g}, \path{u34-poly-group.g}\\[2pt]
Hurwitz action, cusps, the cover $\mathcal B_3$ (\cref{sec:curves,sec:nu}) & finite groups (GAP), a few minutes & \path{u34-census.g}, \path{u34-outer.g}, \path{u34-cusps.g}, \path{u34-marked.g}, \path{u34-classfields.g}\\[2pt]
\Cref{thm:u34}(i), (ii) & exact (PARI/GP); certified continuation of $39$ and $78$ points, $17211$ steps; ball arithmetic for the automorphisms & \path{u34-compute/nu.gp}, \path{u34-nu-verify.py}, \path{u34-aut.py}\\[2pt]
\Cref{thm:u34}(iii) & exact; rank of an elliptic curve (PARI/GP) & \path{u34-compute/E5.gp}, \path{u34-b3-fields.gp}\\[2pt]
\Cref{thm:u34}(iv) & certified continuation of $156$ points, $8907$ steps, $21$ minutes; exact & \path{u34-compute/u34-b3-verify.py}, \path{u34-b3-fields.gp}\\[2pt]
\Cref{thm:poly}(i) & exact (PARI/GP), one second & \path{u34-poly-verify.gp}\\[2pt]
\Cref{prop:monodromy} & certified continuation of $65$ points, $9376$ steps, $47$ minutes; finite groups (GAP) & \path{u34-poly-monodromy.py}, \path{u34-poly-group.g}\\[2pt]
\Cref{prop:arithmetic} & $3179$ factorizations modulo $54$ primes, seconds & \path{u34-poly-check.py}\\[2pt]
\Cref{cor:Q} & exact (PARI/GP) & \path{u34-poly-Q.gp}\\[2pt]
\Cref{prop:rigid,prop:rigidreal} & finite groups (GAP), ten minutes & \path{u34-rigidity.g}\\[2pt]
\Cref{prop:P0} & exact (PARI/GP) & \path{u34-poly-red.gp}\\[2pt]
\Cref{sec:found} (not used in proofs) & numerical & \path{u34-compute}, \path{u34-poly-compute}\\
\bottomrule
\end{tabular}
\makeatletter\def\@captype{table}\makeatother
\caption{The computations.}\label{tab:files}
\end{center}

\subsection{A polynomial of degree \texorpdfstring{$65$}{65} over \texorpdfstring{$\Q$}{Q}}\label{sec:appendix}

The polynomial $P$ of \cref{cor:Q} is recorded in the file \path{u34-poly-Q.txt}; it is not monic, and its coefficients have up to $194$ digits. The function \texttt{polredbest} of PARI/GP \cite{PARI} produces a monic polynomial with smaller coefficients which defines the same number field. Let
\[
P_0(X)=X^{65}+\sum_{k=0}^{64}c_k\,X^k\ \in\Z[X]
\]
be the polynomial whose coefficients $c_k$ are those of \cref{tab:P0}. The largest of them has $127$ digits. The polynomial is also in the file \path{u34-poly-red.txt}.

\begin{proposition}\label{prop:P0}
The polynomial $P_0$ is irreducible, the fields $\Q[X]/(P_0)$ and $\Q[X]/(P)$ are isomorphic, and the Galois group of $P_0$ over $\Q$ is $U_3(4)$, in its action on $65$ points.
\end{proposition}

\begin{proof}
There is a polynomial $a\in\Q[X]$ of degree $64$, recorded in \path{u34-poly-red-map.txt}, with $P(a(X))\equiv0$ modulo $P_0$. This congruence and the irreducibility of $P_0$ are verified in exact arithmetic (\path{u34-poly-red.gp}). Hence the class of $a$ is a root of $P$ in the field $\Q[X]/(P_0)$, and since $P$ is irreducible of the same degree, the two fields are isomorphic. So $P_0$ and $P$ have the same splitting field, and the permutation groups on their roots are isomorphic. The claim follows from \cref{cor:Q}.
\end{proof}

\Needspace{14\baselineskip}
\begingroup\scriptsize
\setlength{\tabcolsep}{2pt}\renewcommand{\arraystretch}{1.0}
\begin{longtable}{@{}r@{\quad}r@{}l@{}}
\caption{The coefficients $c_k$ of $P_0$. Numbers of more than twelve digits are printed in groups of ten digits, from the left, and continue on the following lines.}\label{tab:P0}\\
\toprule
$k$ & & $c_k$\\
\midrule
\endfirsthead
\toprule
$k$ & & $c_k$\\
\midrule
\endhead
\bottomrule
\endfoot
$64$ & $-$ & $19$\\
$63$ & $-$ & $86569$\\
$62$ & $-$ & $4713197$\\
$61$ & $$ & $3870347224$\\
$60$ & $$ & $41455389952$\\
$59$ & $-$ & $3718532384\,3076$\\
$58$ & $$ & $1249182541\,928212$\\
$57$ & $$ & $5210146432\,72600072$\\
$56$ & $-$ & $5454890142\,8671035436\,0$\\
$55$ & $$ & $6185797283\,8120965912\,608$\\
$54$ & $$ & $4379932484\,2572028272\,51632$\\
$53$ & $-$ & $4893980595\,7831960005\,6545484$\\
$52$ & $-$ & $2506300840\,7446786820\,5987736620$\\
$51$ & $$ & $7475056357\,8278924555\,0813239477\,64$\\
$50$ & $-$ & $7359631312\,2153001659\,4732570285\,2396$\\
$49$ & $$ & $2893757103\,4339124177\,0822062876\,210898$\\
$48$ & $-$ & $5036852531\,7528144885\,7888188872\,36181958$\\
$47$ & $$ & $7777061196\,4641106403\,4188436561\,4395599030$\\
$46$ & $-$ & $3364792507\,8287595865\,2657328307\,0226207381\,78$\\
$45$ & $-$ & $3346746030\,2761795134\,0228400612\,0933376213\,0940$\\
$44$ & $$ & $5112634547\,6931634456\,5348331943\,3093230894\,40308$\\
$43$ & $$ & $1158565235\,7986280043\,0078238352\,5870468039\,078923148$\\
$42$ & $-$ & $7602190416\,9463252539\,4177178497\,8704831875\,4178309988$\\
$41$ & $$ & $7153040033\,8242012435\,7349893375\,8142250408\,3063269496\,16$\\
$40$ & $-$ & $2453142749\,2894938011\,7440297931\,6432890520\,6883938452\,68864$\\
$39$ & $$ & $3009550170\,2064770857\,6250773627\,8174028091\,8494124268\,7207320$\\
$38$ & $-$ & $3197294530\,0704926896\,5751473602\,8132759163\,0143220675\,742426616$\\
$37$ & $$ & $2609209549\,9717273036\,5000730750\,1394662118\,7122616756\,8875005484\,4$\\
$36$ & $-$ & $2081290151\,2671772235\,5607126675\,8439837241\,9323955757\,0431833325\,924$\\
$35$ & $$ & $2107276961\,0985144227\,4091369975\,4766177275\,1351241799\,3769946404\,26140$\\
$34$ & $$ & $1254434271\,8095795897\,8348151298\,6103620493\,6239144783\,9255000819\,424764$\\
$33$ & $-$ & $4923384755\,7963579494\,1234490131\,8314347773\,0918586557\,2859346325\,58960443$\\
$32$ & $$ & $1104412218\,4464537140\,8138141288\,6602658034\,1373678972\,7264700335\,7097676241$\\
$31$ & $-$ & $5706685538\,0094601239\,2821911992\,6556623134\,4922365951\,7473781445\,8249485790\,13$\\
$30$ & $-$ & $9727189186\,1391329852\,4956376525\,5350728030\,7953050910\,7403552375\,9141177669\,3353$\\
$29$ & $-$ & $3073568464\,0099388083\,7800105299\,2154162976\,6926071812\,1040614294\,1193056843\,245964$\\
$28$ & $-$ & $2235241251\,3558910392\,9780729839\,6522017175\,2706622699\,0231276344\,6496208185\,19173412$\\
$27$ & $$ & $9865645667\,5400800564\,2149210339\,1552385329\,3333999244\,5464534374\,7834680257\,2167345640$\\
$26$ & $$ & $1469722911\,4085923957\,8111084888\,7713375369\,9468057969$\\*
 & & $7091675722\,2489671423\,0072076314\,624$\\[1pt]
$25$ & $$ & $9759960629\,5845632033\,9780626528\,5975901270\,9600063278$\\*
 & & $4046319215\,8097053244\,7464196973\,6944$\\[1pt]
$24$ & $-$ & $3962113004\,5789178400\,8776428845\,3941229221\,8129223038$\\*
 & & $5965556710\,7510761392\,4431504557\,439552$\\[1pt]
$23$ & $-$ & $1233076942\,2816656108\,7699503095\,8044166566\,1188744360$\\*
 & & $1439878142\,7721883968\,2150932995\,742686720$\\[1pt]
$22$ & $-$ & $8448197832\,1957499294\,8559807215\,2926544223\,3668152616$\\*
 & & $6359069399\,6298532779\,0358224476\,7923391072$\\[1pt]
$21$ & $-$ & $1173020232\,6557855685\,1389024766\,7426513230\,5232415446$\\*
 & & $5543167718\,5170802137\,2766763287\,1049387129\,6$\\[1pt]
$20$ & $$ & $3165973984\,3837429926\,7520636550\,6424294323\,1345548763$\\*
 & & $3071919342\,5267019977\,8288814120\,5196912793\,8048$\\[1pt]
$19$ & $$ & $2074670141\,2233050512\,4390251214\,8892543126\,7670033855$\\*
 & & $7219872979\,8782637343\,8416453001\,7662309934\,677632$\\[1pt]
$18$ & $$ & $2798119518\,3516680527\,4057704482\,4469214802\,7120584844$\\*
 & & $5016459174\,1880677939\,6435488195\,1824752547\,7838080$\\[1pt]
$17$ & $-$ & $3140566296\,5362132608\,2236098279\,0144951952\,3194372168$\\*
 & & $1275371052\,1606444520\,8232463717\,8270376475\,482765312$\\[1pt]
$16$ & $-$ & $2128435689\,4111415312\,3195873579\,8851880199\,1472008605\,6757865776$\\*
 & & $7538077544\,6824616994\,3799258430\,2102195430\,4$\\[1pt]
$15$ & $-$ & $3844318684\,0168154332\,9083997878\,0302081788\,1679409842\,1375190040$\\*
 & & $0173585044\,0706809305\,5772409494\,4936047275\,52$\\[1pt]
$14$ & $$ & $3340262270\,7012120135\,4144822737\,8234656437\,7616239694\,5840264096$\\*
 & & $5733567479\,2968156605\,7536343704\,4509360708\,1728$\\[1pt]
$13$ & $$ & $7055645399\,7603181217\,4373559359\,5612297121\,4544400675\,4205254461$\\*
 & & $5066838093\,9979125951\,5972795059\,8789659261\,84960$\\[1pt]
$12$ & $-$ & $1744423728\,8756067706\,1318070119\,4837732753\,6218169874\,9885830508$\\*
 & & $8181813249\,0462990595\,3637578971\,5068082203\,68124928$\\[1pt]
$11$ & $-$ & $2458520203\,5832314568\,3318292083\,4714549693\,0345886583\,7332433608$\\*
 & & $5498555050\,3539919037\,2145503664\,1846872720\,568995840$\\[1pt]
$10$ & $$ & $3348169058\,0573414187\,8253597057\,9737098675\,4077637059\,5792317979$\\*
 & & $9503980702\,8742384769\,9410786660\,6981978913\,8383479193\,6$\\[1pt]
$9$ & $-$ & $4021052227\,9992374027\,4985212110\,6750993269\,9968083273\,9434759312$\\*
 & & $0898599365\,9571594045\,1571580504\,2197944192\,1631289344\,0$\\[1pt]
$8$ & $-$ & $2184042189\,0268774824\,4095465661\,7006432873\,3172362682\,5397550846$\\*
 & & $3780624296\,3029780718\,1588268544\,8317013968\,7544718329\,4464$\\[1pt]
$7$ & $$ & $1946432688\,3913179230\,4373400006\,8205882894\,3068279841\,4754063081$\\*
 & & $6471774297\,7931068794\,5341068941\,6810345234\,6006730709\,532672$\\[1pt]
$6$ & $$ & $9946446941\,5609572004\,9075339878\,9143301757\,7633048057\,0415201857$\\*
 & & $7503470918\,2819923238\,3937157632\,2572848934\,7638871254\,1585408$\\[1pt]
$5$ & $$ & $4043333033\,7664912102\,5002636501\,5808898996\,8763179848\,1005525644$\\*
 & & $2895962179\,4903980138\,5283909645\,3293490625\,4145538489\,686425600$\\[1pt]
$4$ & $-$ & $1301251670\,6358944031\,7796505137\,4263329278\,3721186386\,7357351964\,7153329471$\\*
 & & $2878530609\,8317635132\,3149305373\,7077744802\,5443690086\,4$\\[1pt]
$3$ & $-$ & $1051593853\,7635656490\,1979723094\,6465040310\,2575172503\,0380813821\,6887787713$\\*
 & & $7968056863\,8324866529\,0544189793\,5280310673\,4189488111\,616$\\[1pt]
$2$ & $$ & $8415731481\,0172693359\,3709568741\,0877487382\,8599825229\,5340357925\,4860553854$\\*
 & & $9405102514\,3001373690\,5152762025\,2594865988\,2913936965\,632$\\[1pt]
$1$ & $-$ & $1926319440\,1861895393\,2512683432\,2329891247\,4824678845\,7679889144\,5793471781$\\*
 & & $1232939290\,8510024250\,7253275764\,3845625000\,9767603171\,49184$\\[1pt]
$0$ & $-$ & $4590892708\,7468796149\,0031076095\,1606098195\,5290843021\,8151272591\,4820061752$\\*
 & & $0462971791\,1188387997\,5109541958\,2265747247\,4213594277\,9764736$\\[1pt]
\end{longtable}
\endgroup

\bibliographystyle{amsplain}
\bibliography{sh-163}

\providecommand{\bysame}{\leavevmode\hbox to3em{\hrulefill}\thinspace}
\providecommand{\MR}{\relax\ifhmode\unskip\space\fi MR }
\providecommand{\MRhref}[2]{%
  \href{http://www.ams.org/mathscinet-getitem?mr=#1}{#2}
}
\providecommand{\href}[2]{#2}
\begin{thebibliography}{10}

\bibitem{BaileyFried2002}
Paul Bailey and Michael~D. Fried, \emph{Hurwitz monodromy, spin separation and
  higher levels of a modular tower}, Arithmetic Fundamental Groups and
  Noncommutative Algebra (Berkeley, CA, 1999), Proc. Sympos. Pure Math.,
  vol.~70, Amer. Math. Soc., Providence, RI, 2002, pp.~79--220.

\bibitem{Atlas}
J.~H. Conway, R.~T. Curtis, S.~P. Norton, R.~A. Parker, and R.~A. Wilson,
  \emph{Atlas of finite groups}, Oxford University Press, Eynsham, 1985.

\bibitem{DebesFried1994}
Pierre D{\`e}bes and Michael~D. Fried, \emph{Nonrigid constructions in {G}alois
  theory}, Pacific J. Math. \textbf{163} (1994), no.~1, 81--122.

\bibitem{FV}
Michael~D. Fried and Helmut V{\"o}lklein, \emph{The inverse {G}alois problem
  and rational points on moduli spaces}, Math. Ann. \textbf{290} (1991),
  771--800.

\bibitem{Johansson2017}
Fredrik Johansson, \emph{Arb: efficient arbitrary-precision midpoint-radius
  interval arithmetic}, IEEE Trans. Comput. \textbf{66} (2017), no.~8,
  1281--1292.

\bibitem{KoenigMC2017}
Joachim K{\"o}nig, \emph{Computation of {H}urwitz spaces and new explicit
  polynomials for almost simple {G}alois groups}, Math. Comp. \textbf{86}
  (2017), no.~305, 1473--1498.

\bibitem{MalleMatzat1999}
Gunter Malle and B.~Heinrich Matzat, \emph{Inverse {G}alois theory}, Springer
  Monographs in Mathematics, Springer, Berlin, 1999.

\bibitem{ShaskaSuperell}
T.~Shaska, \emph{Superellipticity of moduli spaces for covers of the {R}iemann
  sphere}, Preprint, 2026.

\bibitem{GAP}
{The GAP Group}, \emph{{GAP -- Groups, Algorithms, and Programming, Version
  4.12.1}}, 2022, \url{https://www.gap-system.org}.

\bibitem{PARI}
{The PARI~Group}, \emph{{PARI/GP} version \texttt{2.15.4}}, Univ. Bordeaux,
  2023, available from \url{https://pari.math.u-bordeaux.fr/}.

\bibitem{Voelklein}
Helmut V{\"o}lklein, \emph{Groups as {G}alois groups: An introduction},
  Cambridge Studies in Advanced Mathematics, vol.~53, Cambridge University
  Press, Cambridge, 1996.

\bibitem{Wewers1999}
Stefan Wewers, \emph{Deformation of tame admissible covers of curves}, Aspects
  of {G}alois Theory ({G}ainesville, {FL}, 1996), London Math. Soc. Lecture
  Note Ser., vol. 256, Cambridge Univ. Press, Cambridge, 1999, pp.~239--282.

\bibitem{Zywina2025}
David Zywina, \emph{Inverse {G}alois problem for small simple groups}, 2025,
  Note, version of August 2025,
  \url{https://davidzywina.github.io/papers/smallGalois.pdf}.

\end{thebibliography}
\end{document}